\documentclass[12pt]{article}
\usepackage{amsfonts,bbm}
\usepackage{amsthm,comment}
\usepackage{amscd}
\usepackage{amssymb}
\usepackage{hyperref}
\usepackage{amsmath}
\usepackage{tikz} 
\usepackage{etoolbox, textcomp}
\usepackage[matrix,arrow,curve,frame]{xy}

\let\OLDthebibliography\thebibliography
\renewcommand\thebibliography[1]{
  \OLDthebibliography{#1}
  \setlength{\parskip}{0pt}
  \setlength{\itemsep}{0pt plus 0.3ex}
}

\usetikzlibrary{arrows,shapes}
\newcommand{\Uq}{\overline{\text{U}}_q(\mathfrak{sl}_2)}

\newcommand{\Id}{\text{Id}}
\newcommand{\tr}{\text{tr}}

\newcommand{\coend}{\text{coend}\left(\cC\right)}

\newcommand{\bbH}{\mathbb{H}}

\newcommand{\bbC}{\mathbb{C}}   

\newcommand{\bbP}{\mathbb{P}}
\newcommand{\bbZ}{\mathbb{Z}}
\newcommand{\bbX}{\mathbb{X}}

\newcommand{\cP}{\mathcal{P}} \newcommand{\cF}{\mathcal{F}}

\newcommand{\cD}{\mathcal{D}}

\newcommand{\cH}{\mathcal{H}}

\newcommand{\cY}{\mathcal{Y}}  \newcommand{\cZ}{\mathcal{Z}}
\def\cW{\mathcal{W}}
\def\cV{\mathcal{V}}
\def\cC{\mathcal{C}}
\def\cA{\mathcal{A}}  
\def\cN{\mathcal{N}}
  \def\cS{\mathcal{S}}
\def\one{\mathbf{1}_\cC}

\def\Ii{\mathrm{i}}

\newtheorem{theorem}{\textbf{Theorem}}
\numberwithin{theorem}{section}
\newtheorem{lemma}[theorem]{\textbf{Lemma}}
\newtheorem{defn}[theorem]{Definition}
\newtheorem{proposition}[theorem]{\textbf{Proposition}}
\newtheorem{corollary}[theorem]{\textbf{Corollary}}
\newtheorem{remark}[theorem]{Remark}

\newtheorem{conjecture}{\textbf{Conjecture}}
\numberwithin{conjecture}{section}

\theoremstyle{definition}
\newtheorem{ex}{\textbf{Example}}

\newcommand{\hopflink}{{\,\text{\textmarried}}}

\title{Logarithmic conformal field theory, log-modular\\ tensor categories and modular forms\thanks{This review is based on talks we have given in Canberra (July, 2015), Amsterdam (October, 2015), BIRS (February, 
2016) and  Stony Brook (March, 2016).  A more complete draft of our work, including proofs,  is available \cite{CG}.}}

\author{Thomas Creutzig and Terry Gannon}

\date{}

\begin{document}

\tikzset{ar/.style={<-}}
\tikzset{br/.style={->}}

\maketitle

\abstract{
The two pillars of rational conformal field theory and rational vertex operator algebras are modularity of characters on the one hand and its interpretation of modules as objects in a modular tensor category on the other one. Overarching these pillars is the Verlinde formula.

In this paper we consider the more general class of logarithmic conformal field theories and $C_2$-cofinite  vertex operator algebras.
We suggest that their modular pillar are  trace functions with insertions corresponding to intertwiners of the projective cover of the vacuum, and that the categorical pillar are finite tensor categories $\cC$ which are ribbon and whose double is isomorphic to the Deligne product $\cC\otimes \cC^{opp}$.
 Overarching these pillars is then a logarithmic variant of Verlinde's formula. Numerical data realizing this are the modular $S$-matrix and modified traces of open Hopf links.
 
The representation categories of $C_2$-cofinite and logarithmic conformal field theories that are fairly well understood are those of the $\cW_p$-triplet algebras and the symplectic fermions. 
We illustrate the ideas in these examples and especially make the relation between logarithmic Hopf links and modular transformations explicit. 
}

\newpage

\setcounter{tocdepth}{3}

\tableofcontents

\newpage

\section{Introduction}

Vertex operator algebras 
have had a profound influence on mathematics and mathematical physics since their
inception in the mid 1980s. We mention only their roles underlying Moonshine \cite{FLM, Bo},
clarifying conformal field theories,  the chiral de Rham complex in geometry \cite{MSV}, geometrical
Langlands \cite{F}, amongst many others. 

We have a rather complete theory of those vertex operator algebras (VOAs) with 
semi-simple representation theory --- we'll call these \textit{strongly rational} VOAs.
We review this theory in Section \ref{sec:rational} below. Two pillars of this theory are on the one hand the modularity of
characters and, more generally, the automorphicity of chiral blocks for all surfaces, and on the
other hand the interpretation of the modules of the VOA as the objects in a modular tensor category.
Overarching these pillars is the Verlinde formula, which expresses the tensor (fusion) product
coefficients and dimensions of those spaces of chiral blocks, in terms of a quantity (the $S$
matrix) which is common to both pillars.
This picture is completed with a rich collection of examples of strongly rational VOAs.
There are fundamental questions which remain unanswered (e.g. the possible rationality of orbifolds and cosets,
and the possible equivalence with completely rational conformal nets), but over all one must be satisfied
with the picture which has arisen.

We have been working towards formulating a similar picture for more general VOAs. The natural
class to consider is that of $C_2$-cofinite VOAs, which we'll call strongly-finite.  In Section \ref{sec:logVOA} we review their definition and what 
is currently known about them. They correspond to logarithmic conformal field theories.
An analogy to keep in mind is: strongly-rational VOAs are to finite groups, as strongly-finite VOAs
are to finite-dimensional Hopf algebras. We review the representation theory of
Hopf algebras in Section \ref{sec:ass}, where we also explain why representation theory is naturally categorical, and give an abstract meaning
of character. This isn't that dishonest an analogy: to any such VOA we can associate a quantum
group (more precisely, a weak quasi-triangular quasi-Hopf algebra) with identical representation theory.
With this analogy in mind, the theory of strongly-finite VOAs should be a fairly gentle but far-reaching \textit{nonsemi-simple generalization} of 
that of strongly-rational VOAs.
 Nevertheless,  very few nonrational examples have been identified,
and most of the work in the area has focused on studying those isolated examples.
We are trying to see the forest via the individual trees.

The main purpose of this paper is to review the general picture of strongly-finite VOAs 
that is evolving. Section \ref{sec:background} supplies some general background. Section \ref{sec:results} is the heart of
the paper, where we describe our picture. 

For the strongly-rational VOAs, the fundamental general result on modularity is due to Zhu \cite{Z},
who showed that the characters are closed under modular transformations.
Unfortunately this fails in general for strongly-finite VOAs. However, Miyamoto \cite{M1} modified
Zhu's arguments and found that modularity is restored if the characters are supplemented
by finitely many \textit{pseudo-characters}. This is a fundamental result in the theory, but
difficult for many in the theory to understand, and we review it in Section \ref{sec:Miy}. Though
mathematically elegant and natural from the associative algebra perspective, it is not at all
natural from the CFT perspective. Quantum field theory teaches that correlation functions 
should arise from field insertions and ordinary traces. Moreover, it is hard to do computations
within Miyamoto's picture --- it is difficult even to determine in practice the dimension of the resulting
modular group representation. Is there another way to recover Miyamoto's functions, in a manner
more compelling from a standard VOA or CFT point-of-view? We propose in Section \ref{sec:char} a way to do exactly this.

For strongly-rational VOAs, the fundamental result on the second pillar was Huang's proof that
the representation category is a modular tensor category.
The categorical formulation of strongly-finite representation theory, in other words the definition
of a nonsemi-simple analogue of modular tensor category, is more straightforward. We propose
one formulation in Section \ref{sec:logmod}, but another, probably equivalent, is made in
\cite{GaR2}. Thanks primarily to work of Huang and collaborators (see \cite{HLZ8}), the main thing left to prove here is rigidity.
This categorical formulation would immediately imply several interesting things, some of which we
describe in Section \ref{sec:logmod}. 

In the strongly-rational setting, Moore--Seiberg \cite{MS}, Turaev \cite{T}, and Bakalov--Kirillov \cite{BK} explain how to
 obtain from the modular tensor category, a tower of mapping class group representations
and (abstract) spaces of chiral blocks associated to the surfaces with marked points (called
the \textit{modular functor}), and the associated 3-dimensional topological field theory.
This is just the chiral theory; Fuchs, Runkel, Schweigert and collaborators have explained
how to obtain a full conformal field theory from this chiral data \cite{FRS, FFRS}. The analogue of Moore--Seiberg
et al 
for the strongly-finite theory should be Lyubashenko \cite{L,L2,KL}, which we review in Section \ref{lyb},
and Fuchs et al are currently extending their work to the strongly-finite setting \cite{FSS, FS2}.
It is an important question here to show Lyubashenko's representation of SL$(2,\bbZ)$ is
equivalent to that on Miyamoto's space of characters and pseudo-characters.

In the strongly-rational setting, a key to understanding Verlinde (the arch connecting the modular pillar to
the categorical one) is the (numerical) link invariant 
\begin{equation}\label{hopflink}
\begin{tikzpicture}

\draw[ar] (0, 3) arc (40:390:.7cm);
\node[text width=4.2cm] at (-1.8, 2.2)    {$S^\hopflink_{ij}:=\qquad M_i$};

\draw[br] (0, 2) arc (220:570:.7cm);
\node[text width=3cm] at (2.9, 2.3)    {$M_j\quad\in\bbC$};

\end{tikzpicture}
\end{equation}
associated to the Hopf link. In the strongly-finite setting,
too many of these now vanish, and more significant, the fusion rules are no longer diagonalizable in
general so Verlinde's formula must be modified. In the strongly-finite setting, two rings
must be distinguished: the Grothendieck ring spanned by the simple modules, which sees only the
composition factors of all the modules, and the tensor
ring which is spanned by the indecomposable modules and which sees the full splendor of
the tensor product.
We propose in Section \ref{sec:hopfverl} to replace the Hopf link by the \textit{open} Hopf link 
--- the associated link invariants
\begin{equation}\label{openHopf}
\begin{tikzpicture}[baseline=(current  bounding  box.center)]

\draw [ar](.78,1.5) -- (.78,2.1);
\draw [ar](.7, 3) arc (100:435:.4cm);
\draw (.78,2.3) -- (.78,3.5);

\node[text width=3cm] at (1.3, 2.5)    {$U$};
\node[text width=3cm] at (2.1, 1.25)    {$W$};
\node[text width=3cm] at (-1, 2.5)    {$\Phi_{U, W}\ \ \ =$};
\node[text width=3cm] at (3.1, 2.5)    {$\in\mathrm{End}(W)$};

\end{tikzpicture}
\end{equation}
are now matrix-valued and give representations of the tensor ring. These seem sufficient to recover the Jordan blocks in the Grothendieck ring (i.e. the indecomposables in the regular representation 
of the Grothendieck ring) --- this is a purely categorical statement, and should
be possible to prove as in \cite{T}.
To recover the Jordan blocks in the tensor ring, the open Hopf links are not sufficient, but perhaps
knotted versions of the open Hopf link will work. 

In the strongly-finite world, the modular $S$-matrix, which describes how the characters $\chi_M(\tau)$ 
 transform under $\tau\mapsto -1/\tau$, coincides up to a scalar factor with the matrix of Hopf link
 invariants. 
Much more delicate is to see what this relation becomes in the strongly-finite world. Our proposal,
again involving the open Hopf links, is given in Section \ref{sec:conj}.

In Section \ref{sec:V} we show what our conjectures look like for  the $\cW_p$ models
and in Section \ref{SF} we discuss symplectic fermions $SF_d^+$.

Most urgently, we need more independent examples of strongly-finite VOAs. Perhaps this limits the development of the theory more than anything else.
We address some possibilities for this in Section \ref{sec:ex}. Again, there should be a far richer zoo
of strongly-finite VOAs than strongly-rational ones, but at present the opposite is emphatically
the case.

We are aware of the reluctance, not only by physicists but by many mathematicians, towards categorical formulations.
In fact we have little sympathy for abstract nonsense done for its own sake. But in representation
theory especially, it seems an indispensable tool: it helps obtain interesting results that don't
need category theory to formulate. In our concluding section we supply some
evidence to the skeptic that category theory is a practical tool for VOAs and CFT.

\vspace{5mm}

\noindent {\bf Acknowledgements:}
We thank Azat Gainutdinov and Ingo Runkel for discussions on related topics and for sharing \cite{GaR2} which has partial overlap with our results.
 We would like to thank  J\"urgen Fuchs for many explanations of Lyubashenko's work on modularity within nonsemi-simple 
tensor categories. 
 TC also thanks Shashank Kanade for many discussions on related topics.  TG thanks the Physics Department at Karlstad University for a very pleasant and
stimulating work environment. Our research is supported in part by NSERC.

\section{Background}\label{sec:background}

\subsection{Representation theory of  associative algebras}\label{sec:ass}

We will assume the reader is familiar with  the finite-dimensional representations over $\bbC$ of a finite group $G$. The theory is \textit{semi-simple}:
$G$ has finitely many irreps 
$\rho_1,\ldots,\rho_r$ (up to equivalence), and every finite-dimensional representation $\rho$ is (up to equivalence) the direct sum $\rho\cong
\oplus_{i=1}^rm_i\rho_i$ of irreps in a unique
way. We have a tensor product $\rho\otimes\rho'$ of representations, which (up to equivalence at least) we can
completely capture by the tensor product multiplicities $T_{ij}^k$ defined by $\rho_i\otimes\rho_j\cong\oplus_kT_{ij}^k\rho_k$. 
Letting $[\rho]$ denote the equivalence class of $\rho$, and writing $[\rho\oplus\rho']=[\rho]+[\rho']$ and $[\rho\otimes\rho']
=[\rho][\rho']$, we get a ring structure on the $\bbZ$-span of the irreps $[\rho_i]$, with  structure constants 
 $T_{ij}^k$. We can call this ring the \textit{Grothendieck ring} of $G$. We  also have a  trivial 1-dimensional representation $\rho_1={\bf 1}$, which is the tensor unit.
And we have duals $\rho^*$, which we can completely capture (up to equivalence) by an order 1 or order 2 permutation $i\mapsto i^*$ on the index set $\{1,\ldots,r\}$. All of this information is no more and no less than the character table
of $G$ itself.

The usefulness of character tables to the representation theory of finite groups is clear. But different groups (e.g. the dihedral group $D_4=\langle a,b\,|\,a^2=b^4=abab=1\rangle$ and the group $Q_8=\langle i,j\,|\,i^2=j^2,\,i^4=1,\,iji=j\rangle$)
can have the same character table. How can we enhance
the character theory, to get something representation theoretic which does a better job characterizing $G$?

It should be obvious how we can enhance this data. The key phrase, used over and over two paragraphs ago, was `up to equivalence'. That was why it reduced to the character table. 
But we don't just have characters, we also have {\it intertwiners} $T\in\mathrm{Hom}_G(\rho,\rho')$, i.e. linear maps $T:\bbC^{d'}
\rightarrow\bbC^d$ (where $d,d'$ are the dimensions of $\rho,\rho'$ respectively) satisfying $\rho\circ T=T\circ\rho'$. Put another way, what we really have is a {category}: the
objects are the  representations (or modules) and the arrows or morphisms between the objects are the intertwiners.
This category has lots of additional structure: it has direct sums, tensor products, complete reducibility, duals, etc.
The result is called a tensor category (in fact a fusion category, which is even better). We  review basic categorical notions in Section \ref{sec:btc}.

In fact, if we think about it, character tables can be surprisingly ineffective. The representation theory of the groups $D_4$ (dihedral) and $Q_8$ (quaternions)  actually have many differences. For one thing, it is elementary that
the determinant of a $G$-representation, i.e. $(\mathrm{det}\,\rho)(g)=\mathrm{det}\,\rho(g)$ will always be a 1-dimensional $G$-representation. The determinant of the unique 2-dimensional irrep of $D_4$ is nontrivial, whereas the determinant of the 
unique 2-dimensional irrep of $Q_8$ is trivial. So  $D_4$ and $Q_8$ could be 
distinguished if we also knew the determinants of the irreps. We've just learned that character  tables can't see the determinant of
irreps! Also, the 2-dimensional irrep of $D_4$ can be realised as
\begin{equation}\label{rep1}a\mapsto \left(\begin{matrix}1&0\cr 0&-1\end{matrix}\right)\,,\ b\mapsto \left(\begin{matrix}0&-1\cr 1&0\end{matrix}\right)\end{equation}
that is, as a real orthogonal matrix representation. Such a representation {\it equals} its own dual, not merely is equivalent to
its dual. It is easy to see that this cannot be done for $Q_8$. So the character table by itself cannot determine if representations
are self-dual, but only whether they are equivalent to their duals.

But can these things be seen by tensor categories? Another realisation of the 2-dimensional $D_4$-irrep   is
$$a\mapsto \left(\begin{matrix}0&1\cr 1&0\end{matrix}\right)\,,\ b\mapsto \left(\begin{matrix}i&0\cr 0&-i
\end{matrix}\right)$$
It satisfies $\rho^*\cong\rho$, but in the tensor category we would also know the intertwiner. In fact in the
tensor category the equivalence $\rho^*\cong \rho$ is replaced with the {\it equality} $\rho^*=\psi\otimes\rho$,
where $\psi$ is the 1-dimensional representation defined by $\psi(a)=1,\psi(b)=-1$. However, the matrix representation in \eqref{rep1}
is a different object in the category (even though it is an equivalent representation), and for it we obtain $\rho^*=\rho$.
In both cases the equivalence $\rho\cong\rho^*$ or $[\rho]=[\rho^*]$ is replaced by something stronger in the category.
Tensor categories can also see the determinant det$\,\rho$, using exterior power. 

 So the tensor categories can distinguish $D_4$ and $Q_8$.
However the tensor category alone does not distinguish all groups: e.g. there are groups of order $2^{15}3^{4}\cdot5\cdot7\approx92$ million which are identical as tensor categories (surely there are much smaller examples). 

However, there is one more piece of information which can be included. We know $\rho\otimes \rho'\cong \rho'\otimes \rho$
(after all, they have identical characters $\chi_\rho\chi_{\rho'}=\chi_{\rho'}\chi_\rho$). So that means there is an invertible intertwiner
$c_{\rho,\rho'}\in {\rm Hom}_G(\rho\otimes \rho',\rho'\otimes \rho)$. 
 For finite groups,
this is just $c_{\rho,\rho'}(u\otimes v)=v\otimes u$. Of course this obeys $c_{\rho,\rho'}\circ c_{\rho',\rho}=Id_{\rho'\otimes \rho}$.
This defines what is called a symmetric tensor category.
 {\it Tannaka-Krein duality} says that the tensor category of $G$ interpreted as a symmetric tensor category, uniquely determines $G$.

The moral: in representation theory, consider the {\it category} of representations, not just the representations up to equivalence. 
The category of representations of $G$ should be regarded as a gentle enhancement of the combinatorics of 
representations, which carries much more information. In short: \textit{representation theory is categorical.}

The representations of finite groups can be recast as modules of the group algebra $\bbC G$. More generally,
one can consider the finite-dimensional modules over an associative algebra $\cA$. We know we should consider intertwiners
Hom$_\cA(U,V)$, i.e. linear maps $T:V\rightarrow U$ satisfying $a.T(v)=T(a.v)$ for all $a\in\cA$ and $v\in V$, i.e. we
should consider the
category {Mod}$^{fin}(\cA)$ of $\cA$-modules. $\cA$ will have finitely many irreducible
modules, and direct sum of modules will exist, but semi-simplicity (complete reducibility) of modules will usually be lost. If $\cA$ is e.g. a Hopf algebra
(as $\bbC G$ is), then one will be able to define 
 tensor products of modules using comultiplication, a  tensor unit
(i.e. an analogue of  trivial representation) using the counit, and each module will have duals
(using the antipode). In fact  {Mod}$^{fin}(\cA)$ is an example of a \textit{finite tensor category} (discussed in
Section \ref{sec:btc}).  If e.g. the Hopf algebra $\cA$ is cocommutative (as $\bbC G$ is), we should consider the intertwiners $c_{U,V}$
for all $\cA$-modules $U,V$. But when the intertwiners $c_{U,V}$ do exist, they may not satisfy
$c_{U,V}\circ c_{V,U}=Id_{V\otimes U}$. In fact this latter possibility turns out to be {\it by far} the most interesting: {Mod}$^{fin}(\cA)$ in this case is called a
{\it braided tensor category}. We also discuss these in Section \ref{sec:btc}.

We are all familiar with the notion of \textit{character} of a finite group $G$.  They 
 form a basis for
the space of \textit{class functions} $f(h)$, maps $f:G\rightarrow \bbC$ satisfying $f(k^{-1}hk)=f(h)$.
The associative algebra associated to a group, of course,
is  its group algebra $\bbC G=\{\sum_{g\in G}c_gg\,|\,c_g\in\bbC\}$. Then $\bbC G$ is a bimodule over itself, or equivalently it carries a representation of $G\times G$:
it acts on itself  by both left- and  right-multiplication. This bimodule is
$\bbC G^{bi}\cong \oplus_MM\times M^*$, where the sum runs over the finitely many irreducible $G$-modules,
and $M^*$ denotes the dual or contragredient. This bimodule $\bbC G$ is central to the
whole representation theory of $G$: we can restrict this $G\times G$ representation
to $G\times 1$ or $1\times G$, in which case it is the (left- or right-)regular representation
$\bbC G^{reg}\cong \oplus_M\mathrm{dim}\,M\,M$. Or we can restrict it to the diagonal subgroup
$\{(g,g)\}\cong G$ in $G\times G$, in which case we get the adjoint representation $\bbC G^{adj}
\cong \oplus_M M\otimes M^*$.

Abstractly, a class function for the group $G$
is precisely an \textit{intertwiner}, an element $f \in \textrm{Hom}(\bbC G^{adj},\bbC)$. This interpretation of character
(more precisely, of class function) underlies the space in Lyubashenko on which the modular group
is to act, as we'll see in Section 2.7 below.

Most associative algebras $\cA$ will not have a semi-simple representation theory. In this case each $\cA$-module
will decompose in a unique way as a direct sum of \textit{indecomposables}, but not all
indecomposable modules will be simple (in fact generic finite-dimensional $\cA$ will have
uncountably many indecomposable modules but only finitely many simple ones).

For example, the simple modules of the polynomial algebra $\cA=\bbC[x]$ are all 1-dimensional:
$x$ acts like a number $\lambda$ on $\bbC^1$. But it has indecomposable modules in
each dimension, where $x$ acts on $\bbC^d$ by a Jordan block $B_{\lambda,d}$.

Given an $\cA$-module $M$, by its \textit{socle} we mean the largest semi-simple submodule.
The \textit{Loewy diagram} of $M$ has bottom row the socle $soc(M)$ of $M$, on top of that is put the
socle of  $M/soc(M)$, etc. Writing the Loewy diagram horizontally essentially gives the 
composition series of $M$: $0\hookrightarrow M^1\hookrightarrow M^2\hookrightarrow\cdots
\hookrightarrow M^k=M$ where $M^{i}$ are all submodules of $M$ and each $M ^{i+1}/M^i$
is simple (called a \textit{composition factor} of $M$). The composition factors are the union (keeping
multiplicities) of all simple modules appearing in each row of the Loewy diagram.

For example the Loewy diagram of the $x\mapsto B_{\lambda,d}$ 
indecomposable of $\bbC[x]$ consists of a vertical line of length $d$ with the $x\mapsto \lambda$
simple in each spot. The composition factors are $\lambda$ (with multiplicity $d$).

When semi-simplicity is lost, much of the role played by the simple modules in the semi-simple case
is shared by simple modules and their projective covers.  
Given any algebra $\cA$ (e.g. a VOA), an $\cA$-module $P$ is called \textit{projective} if for every $\cA$-modules $M,N$
and surjective intertwiner $f:M\rightarrow N$ and any intertwiner $h:P\rightarrow N$, there is an
intertwiner $h':P\rightarrow M$ such that $h=f\circ h'$.  An  intertwiner $h:P\rightarrow M$ is called a \textit{projective cover} of an
$\cA$-module $M$ if $P$ is a projective $\cA$-module and any intertwiner $g:N\rightarrow P$ for which $h\circ g$ is
surjective, is itself surjective.  

For example, the simple modules of $\bbC[x]$ have no projective cover. Choosing instead
$\cA=\bbC[x]/(x^3)$ say, the only  simple module is $x\mapsto 0$, and its projective cover is
$x\mapsto B_{0,3}$. In this example, up to isomorphism the only other indecomposable is
$x\mapsto B_{0,2}$.

\subsection{VOA basics}\label{sec:basics}

We assume the reader has at least a vague understanding of the definition of a VOA and its
modules. The purpose of this short subsection is to collect together
some more technical aspects which will help the reader through some of the later material.

For reasons of simplicity, we will restrict attention in this paper to  vertex operator algebras (VOAs) $\cV$ of
CFT-type --- this means $\cV$ has $L_0$-grading $\cV=\coprod_{n=0}^\infty \cV(n)$ with $\cV(0)=\bbC\mathbf{1}$
and dim$\,\cV(n)<\infty$. We also require $\cV$, when regarded as a $\cV$-module in the usual way, to be simple and 
isomorphic to its contragredient $\cV^\vee$. For most of the paper we'll require $\cV$ to be 
 $C_2$-cofinite, a standard condition implying $\cV$ has finitely many simple modules,
 a tensor product, convergent characters, etc.

 Unless otherwise stated, we'll assume that $\cV$ is a \textit{strongly-finite} VOA, by which we mean a simple $C_2$-cofinite  VOA of CFT-type with $\cV$ isomorphic to $\cV^\vee$. If $\cV$ is
 in addition regular, then $\cV$ is called \textit{strongly-rational}.

By Mod${}^{g.r}(\cV)$ we mean the category of  \textit{grading-restricted} weak $\cV$-modules $M$,
by which we mean $L_0$ decomposes $M$ into a direct sum of finite-dimensional generalized eigenspaces $M(h)$ with
eigenvalues $h$ bounded from below. The \textit{conformal weight} $h_M$ of $M$ is the 
eigenvalue of smallest real part, when it exists. When $L_0$ is in fact diagonalizable over $M$,  $M$ is called \textit{ordinary}.
For convenience we will refer to any grading-restricted
weak $\cV$-module, simply as a $\cV$-module.

Let $M,M'$ be $\cV$-modules. By a \textit{homomorphism} $f:M\rightarrow M'$ we mean
a linear function satisfying $f(Y^M(v,z)w)=Y^{M'}(v,z)f(w)$ for all $v\in\cV$ and
$w\in M$, or equivalently $f\circ v^M_n=v^{M'}_n\circ f$ for all $v\in \cV$ and $n\in \bbZ$.
 For example, $L_0$ is always
in the centre of End$(M)$.

Let $M,M',M''$ be any  $\cV$-modules. The notion of \textit{logarithmic intertwiner $\cY$ of type $\left({M''\atop M,M'}\right)$} was introduced in \cite{Mil}, building upon \cite{FHL}. Among other things,  for each $u\in M$, $\cY$ has an expansion
 $$\cY(u,z)=\sum_{i=0}^N\sum_{n\in h+h'-h''+\bbZ}u_n^{(i)}z^{n-1}(\mathrm{log}\,z)^i$$
 for some $N$ depending on $\cY$,
 where for $u\in M(k)$, each $u_n^{(i)}$ maps the generalized $L_0$-eigenspace $M'(\ell)$
 into $M''({k+\ell-n-1})$.

Let $M,N$ be any indecomposable $\cV$-modules with conformal weights $h_M,h_N$,
and choose any intertwiner $\cY$
 of type $\left({M\atop N,M}\right)$. Then for each $u\in N$, $\cY$ has an expansion
 $$\cY(u,z)=\sum_{i=0}^K\sum_{n\in h_N+\bbZ}u_n^{(i)}z^{n-1}(\mathrm{log}\,z)^i$$
 where for $u\in N(k)$, each $u_n^{(i)}$ maps the generalized $L_0$-eigenspace $M(\ell)$
 into $M({k+\ell-n-1})$. Write $o_\cY(v)=v_{k-1}^{(0)}$ for $v\in N(k)$. 
For any $v\in N$ and $\tau$ in the upper half-plane $\bbH$, define the \textit{1-point function}
 \begin{equation}Z^\cY(v,\tau)=\mathrm{Tr}_M(o_\cY(v)q^{L_0-c/24})\,.\label{1ptfn}\end{equation}
These are central to the modularity of our story. We restrict here to intertwiners of type 
$\left({M\atop N,M}\right)$, as otherwise the traces wouldn't mean anything.

An important case of \eqref{1ptfn} is $N=\cV$. In this case the intertwiner space is 1-dimensional, with basis given by the map $Y_M$ defining the $\cV$-action on $M$.
Then \begin{equation}\label{char}\mathrm{ch}[M](v,\tau):=Z^{Y_M}(v,\tau)=q^{-c/24}\sum_{n=0}^\infty\mathrm{Tr}_{M(n+h_M)}o(v)\,q^{n+h_M}\,,\end{equation}
where  $o(v)\in \text{End}(M)$ is the zero-mode of the field of $v$,  is called the \textit{character} of $M$. Most important is the specialisation
$v=1$, which we will call the \textit{graded-dimension} of $M$:
\begin{equation}\mathrm{ch}[M](\tau):=Z^{Y_M}(1,\tau)=q^{-c/24}\sum_{n=0}^\infty\mathrm{dim}({M({n+h_M})})\,q^{n+h_M}\label{graded-dim}\end{equation}
The term `character' is usually used for \eqref{graded-dim}, but this is misleading
considering the analogy with finite groups.

 Write $[M]$ for the equivalence class of $\cV$-modules 
isomorphic to $M$. Define $[M]+[N]=[M\oplus N]$, $[M][N]=[M\otimes N]$, and $[M]^*=[M^*]$.
Extend all this in the usual way to formal combinations over the integers, and
the result is the tensor product ring, which we'll denote Fus$^{full}(\cV)$. In the nonrational case,
this is a very big ring, and for most or all purposes it suffices to consider the subring Fus$^{simp}(\cV)$ generated
by both simple and projective modules. In all strongly-rational examples we have seen, this ring is finite-dimensional
(over $\bbZ$), and we expect this to hold in general.

A different smaller version of Fus$^{full}(\cV)$ is the \textit{Grothendieck ring} Fus$^{gr}(\cV)$. 
Recall that each $\cV$-module $M$  has a composition series, which describes
how to build it up by extending by simple modules, called its composition factors. Write $[[M]]$ for the formal sum 
of all equivalence classes of composition factors of $M$ (with multiplicities). Then the definitions $[[M]]+[[N]]:=[[M\oplus N]]$,
$[[M]][[N]]:=[[M\otimes N]]$ are well-defined and define a ring structure on the
$\bbZ$-span of the simple classes $[[M_j]]$. 

 We prefer to avoid calling any of these  \textit{the} fusion ring, as the term
 `fusion ring' is used in different senses in the literature.

\subsection{The rational modular story}\label{sec:rational}

Let $\mathcal V$ be strongly-rational, then such a VOA has the following three properties:

\begin{itemize}

\item[\textbf{(Cat)}] its category of modules Mod$^{g.r}(\cV)$
is a modular tensor category \cite{T,H2};  

\item[\textbf{(Mod)}] its torus 1-point functions \eqref{1ptfn} are modular \cite{Z,M0}; and 

\item[\textbf{(Ver)}] Verlinde's formula
is true \cite{V,MS,H1}.
\end{itemize}

This will now be explained in more detail, starting with the categorical aspect:

\smallskip
\hspace{-1cm} \textbf{(Cat)}: In this case all $M\in\mathrm{Mod}^{g.r}(\cV)$ are ordinary and there are only finitely many isomorphism classes of simple modules.
We choose representatives $M_1, \dots, M_n$ of these classes, with $M_1=\cV$ the VOA itself. Every (ordinary) $\cV$-module will
 be isomorphic to a unique direct sum of these simple $M_j$.  The \textit{fusion coefficients} $\cN_{ij}{}^k\in\bbZ_{\ge 0}$ are nothing but the dimensions of the spaces of intertwining operators of type $\left({M_k\atop M_i\,,M_j} \right)$.
 
Huang and Lepowsky have been able to define tensor products of modules (see e.g. \cite{HL3} and references therein), realizing those fusion rules:
\[M_i\otimes M_j\cong \bigoplus_{k=1}^n {\cN_{ij}}^k\, M_k\,.\]
 This tensor product gives a 
ring structure (the \textit{tensor  ring} or \textit{fusion ring}) on the formal integer-span of the representatives
$M_j$, with unit $M_1$. For each $i$, the fusion matrix is the matrix representing the action of $M_i$ on modules via tensor product:  $\cN_i$ is defined by $(\cN_i)_{jk}={\cN_{ij}}^k$. They yield a
 representation of the tensor  ring:
 \[\cN_i\,\cN_j=\sum_k{\cN_{ij}}^k\,\cN_k\,.\]
Fusion coefficients are also dimensions of spaces of linear maps, $\cN_{ij}{}^k=\mathrm{dim(Hom}_\cV(M_i\otimes M_j,M_k))$.
 
We turn to the modular aspect:

\smallskip
\hspace{-1cm} \textbf{(Mod):} For a $\cV$-module $M$, recall the character ch$[M](v,\tau)$,
defined for any $v\in\cV$ and $\tau$ in  the upper half-plane $\bbH$.  
Zhu has shown \cite{Z} that for fixed $v\in\cV{[k]}$ (using the $L_{[0]}$-grading introduced in \cite{Z}), that these ch$[M](v,\tau)$
 form the components of a vector-valued modular form of weight $k$ for the modular group SL$(2, \bbZ)$. Further they are holomorphic 
in $\bbH$ and meromorphic at the cusps. Especially, the graded-dimensions form a
vector-valued modular form of weight 0 for SL$(2,\bbZ)$.
We stress that in many examples graded-dimensions are not linearly independent but characters are. This is the reason for preferring the more general characters over the more familiar graded-dimensions.
This modularity (and linear independence) means, for any $v\in\cV{[k]}$, the modular $S$-transformation $\tau\mapsto -1/\tau$ defines a matrix 
$S^\chi$ via
\[
\mathrm{ch}[{M_i}](v,-1/\tau) = \tau^k\sum_{j=1}^n  S^\chi_{ij}\,\mathrm{ch}[{M_j}](v,\tau)\,.
\]

More generally, choose any $\cV$-module $N$ and vector $v\in N{[k]}$. Then Miyamoto showed that Zhu's methods naturally generalize to show that 
the 1-point functions $Z^{\cY}(v,\tau)$, as $\cY$ runs over a basis of intertwiners of type $\left({M\atop N,M}\right)$ for all simple modules $M$, form a vector-valued modular form of weight $k\in h_N+\bbZ$
for SL$(2,\bbZ)$. 

The eigenvalues of $L_0$ on a simple $\cV$-module $M_j$ are elements of $h_j+\bbZ_{\ge 0}$ where $h_j$ is the
{conformal weight} of $M_j$. Modularity of characters implies that these $h_j$, together with the central charge $c$,  are rational \cite{AnMo}.

\smallskip
\hspace{-1cm} \textbf{(Ver):}
When $\cV$ is strongly-rational, the Grothendieck ring Fus$^{gr}(\cV)$, full tensor
 ring Fus$^{full}(\cV)$ and simple-projective subring Fus$^{simp}(\cV)$ all coincide.
In 1988, Verlinde \cite{V}  conjectured that these so-called fusion coefficients ${\cN_{ij}}^k$ are related to the 
matrix $S^\chi$ by what is now called \textit{Verlinde's formula:}
\begin{equation}\label{VerlForm}
{\cN_{ij}}^k = \sum\limits_{\ell=1}^n \frac{S^\chi_{i\ell}S^\chi_{j\ell}\left(S^{\chi\,-1}\right)_{\ell k}}{S^\chi_{1\ell}}\,.
\end{equation}
This is one of the most exciting outcomes of the mathematics of rational conformal field theory
(CFT).  
As $S^\chi$ is a unitary matrix, $(S^{\chi\,-1})_{\ell k}$ here can be replaced with the complex conjugate
 $S^{\chi\,*}_{k\ell}$.  The Verlinde formula thus has three aspects:
 
 \begin{itemize}

\item[\textbf{(V1)}] There is a matrix $S^{\otimes}$ simultaneously diagonalizing all fusion matrices $\cN_{i}$. Each diagonal entry $\rho_l(M_i):=(S^{\otimes\,-1}\cN_{i}S^{\otimes})_{ll}$
defines a one-dimensional representation 
\[\rho_l(M_i)\rho_l(M_j)=\sum_k{\cN_{ij}}^k\rho_l(M_k)\]
of the tensor  ring. All of these $\rho_l$ are distinct.

 \item[\textbf{(V2)}] The Hopf link invariants \eqref{hopflink}
 for any $1\le i,j\le n$ give one-dimensional representations of the fusion ring: 
\[\frac{S^\hopflink_{i\ell}}{S^\hopflink_{1\ell}}\frac{S^\hopflink_{j\ell}}{S^\hopflink_{1\ell}}=\sum_{k=1}^n{\cN_{ij}}^k\frac{S^\hopflink_{k\ell}}{S^\hopflink_{1\ell}}\,.\]
Moreover, each representation $\rho_{l}$ appearing in (V1) equals one of these $M_i\mapsto {S^\hopflink_{i\ell}}/{S^\hopflink_{1\ell}}$.

\item[\textbf{(V3)}] \textit{The deepest fact:} these three $S$-matrices are essentially the same. More precisely,
${S^\chi_{ij}}/{S^\chi_{11}}=S^\hopflink_{ij}$, and $S^{\otimes}=S^\chi$ works in (V1). For this choice of
$S^{\otimes}$, $\rho_l(M_i)={S^\hopflink_{il}}/{S^\hopflink_{1l}}$.

\end{itemize}

Much of (V1)-(V3) is automatic in any modular tensor category. In particular, Theorem 4.5.2 in \cite{T} says that in
any such category,
\[\cN_{ij}^k=\cD^{-2}\sum_\ell\frac{S^\hopflink_{i\ell}S^\hopflink_{j\ell}S^{\hopflink\,*}_{k\ell}}{S^\hopflink_{1\ell}}\,,\]
where $\cD^2=\sum_iS^{\hopflink\,2}_{1i}$. So the remaining content of \eqref{VerlForm} is that $S^\chi=\cD^{-1}S^\hopflink$. 
Categorically, the space of 1-point functions $F_M$ is identified with the space Hom$(1,\cH)$ where $$\cH\cong
\bigoplus_iM_i^\vee\otimes M_i$$ (throughout this paper, $M^\vee$ denotes the \textit{contragredient} or dual of 
$M$). This should remind us of the definition of class function from Section 2.1, and that is no accident. This object $\cH$ is naturally a Hopf algebra and a Frobenius algebra. Using this structure,
Section 6 of \cite{L} shows End$(\cH)$  carries a projective action of SL$(2,\bbZ)$; it is the coend of \cite{L}, and \cite{Sh} suggests to interpret it (or its dual) as the adjoint algebra of the category. The SL$(2,\bbZ)$-action
on the $F_M$ correspond to the subrepresentation on  Hom$(1,\cH)$,
which is generated by $S$-matrix $\cD^{-1}S^\hopflink$ and $T$-matrix coming from the ribbon twist. The other subrepresentations Hom$(M_j,\cH)$
correspond to the projective SL$(2,\bbZ)$-representations on the space of 1-point functions with insertions $v\in M_j$, defined in \eqref{1ptfn}. As in \cite{FS}, the character ch$[M]$ of \eqref{char} can be interpreted categorically as
the partial trace of the `adjoint' action $\cH\otimes M\rightarrow M$.

It turns out to be far easier to construct modular tensor categories, than to construct VOAs or CFTs. This is related to the
fact that the subfactor picture is essentially equivalent to the categorical one, at least if one assumes unitarity. This can be 
used to probe just how complete the lists of known (strongly-rational) VOAs and CFTs are. More precisely, given a fusion category
(essentially a modular tensor category without the braiding), taking its Drinfeld double or centre construction
yields a modular tensor category. Fusion categories can be constructed and classified relatively easily, and their doubles
worked out, at least when their fusion rules are relatively simple. Examples of this strategy are provided in e.g. \cite{EG0}.
This body of work suggests that the zoo of known  modular tensor categories is quite incomplete, and hence that the zoo of
known  strongly-rational VOAs and CFTs is likewise incomplete.

\subsection{The finite logarithmic story}\label{sec:logVOA}

An important challenge is to extend the aforementioned results beyond the semi-simplicity of the associated tensor categories. 
A VOA or CFT is called \textit{logarithmic} if at least one of its modules is indecomposable but reducible. The name refers to 
logarithmic singularities appearing in their correlation functions and operator product algebras of intertwining operators. In this
paper we are primarily  interested in strongly-finite VOAs, i.e. logarithmic VOAs with only finitely many simple modules (see Section 2.2 for the formal definition). Of these, the best studied class is the family of $\cW_p$-triplet algebras parameterized by $p\in\mathbb Z_{>1}$ and this will be our main example too. 
See e.g. \cite{CR4,FS} for introductions on logarithmic VOAs.

Though strongly-finite VOAs have finitely many simple modules, they usually have uncountably many indecomposable ones. Any strongly-rational
VOA is strongly-finite. 
This subsection reviews what is known about them. The key hypothesis here
is $C_2$-cofiniteness; many of the following aspects will
persist even when simplicity or CFT-type is lost.

The $\cW_p$ models are strongly-finite. They have central charge $c=1-6(p-1)^2/p$, and are generated by the conformal
vector and 3 other states. The symplectic fermions form a 
logarithmic vertex operator superalgebra with $c=-2d$ for any $d\in\bbZ_{>0}$ (the number of pairs of fermions); their even part $SF^+_d$ is a strongly-finite VOA \cite{Abe}. The $\cW_{p,p'}$-models are $C_2$-cofinite \cite{TW2} but not simple, so are not
strongly-finite.

\begin{itemize}

\item[\textbf{(Cat)$'$}]  
A tensor product theory for strongly-finite VOAs has been developed by Huang, Lepowsky and Zhang (see e.g. \cite{HLZ8} and references therein), see also Miyamoto \cite{M2}. The corresponding tensor category Mod$^{g.r}(\cV)$ is braided. 
Rigidity of this category is   proven so far only for the $\cW_p$ models \cite{TW} and symplectic fermions $SF^+_d$ \cite{DR}.\end{itemize}

There have been different proposals in the literature for what replaces modular tensor category here
(e.g. \cite{KL,FS}. See Section \ref{sec:logmod} for our preference.
 In this nonsemi-simple setting, 
we will in general only have the inequality $\cN_{U,V}{}^W\le \mathrm{dim(Hom}_\cV(U\otimes V,W))$.

Fix representatives $M_1,\ldots,M_n$ of isomorphism classes of simple $\cV$-modules as before.  Let $P_i$ be their projective covers. Likewise
fix representatives of isomorphism classes $M_\lambda$ of indecomposable $\cV$-modules. Any $\cV$-module $M\in\mathrm{Mod}^{g.r}(\cV)$ is
isomorphic to a direct sum of finitely many $M_\lambda$, in a unique way. 
Recall from Section 2.2 the  \textit{(full) tensor  ring} Fus$^{full}(\cV)$, the  subring Fus$^{simp}(\cV)$, and the  the \textit{Grothendieck ring} Fus$^{gr}(\cV)$.
These will no longer be isomorphic in general --- in fact Fus$^{full}(\cV)$ will usually have an uncountable
basis. We prefer to avoid calling any of these  \textit{the} fusion ring, as the term
 `fusion ring' is used in different senses in the literature.

 A serious problem here is that the $\bbC$-span of graded-dimensions ch$[M_i](\tau)=F_{M_i}(\tau,1)$
is   no longer SL$(2,\bbZ)$-invariant, although we still have that each ch$[M_i]\left(\frac{a\tau+b}{c\tau+d}\right)$ lies in
the $\bbC[\tau]$-span of the characters. Many authors interpret this as saying that e.g. $S^\chi$ is now $\tau$-dependent.
We prefer to say that the \textit{ordinary} trace functions (or characters) \eqref{char} must  be augmented by  \textit{pseudo-trace} functions (or pseudo-characters)
associated to indecomposable $\cV$-modules. The main result here is due to Miyamoto \cite{M1}:

\begin{itemize}

\item[\textbf{(Mod)$'$}] The $\bbC$-span of (pseudo-)characters (suitably defined)
 is SL$(2,\bbZ)$-invariant.
\end{itemize}

Corollary 5.10 of \cite{M1} uses this to argue that conformal weights and central charges in strongly-finite VOAs must also be 
rational. The dimension of the resulting SL$(2,\bbZ)$-representation is dim$A_m/[A_m,A_m]-\mathrm{dim}A_{m-1}/[A_{m-1},
A_{m-1}]$, where $A_k=A_k(\cV)$ is the $k$th Zhu algebra (a finite-dimensional associative algebra), and $m$ is sufficiently large. It is possible to state explicitly how large $m$ must be, but this isn't terribly useful as the algebras $A_k(\cV)$ are very hard to identify in practise. We review this
 in Section \ref{sec:Miy}.

The conformal field theory philosophy is that traces with field insertion are the natural objects to look at. This point of view has first appeared in the context of boundary conformal field theory of $bc$-ghost and  then symplectic fermions \cite{CQS, CRo} and then later in the symplectic fermion super VOA \cite{Ru}. Computations of them in the spirit of Miyamoto using symmetric linear functions 
has so far also only be done in the symplectic fermion case \cite{AN}.

There have been  several proposals for a strongly-finite Verlinde formula beyond rationality, mainly in the physics literature, see \cite{R, GR, FHST,GaiT,FK}. All these proposals are guided by analogy to the rational setting and they do not connect to the tensor category point of view (recall that equation \eqref{VerlForm} is a theorem in any modular tensor category, when $S^\chi$
there is replaced with $\cD^{-1}S^\hopflink$). 
Define the tensor  resp. Grothendieck matrices $\cN_M$ resp. $\cN^{gr}_{i}$, with
entries $(\cN_M)_{M',M''}={\cN_{M,M'}}^{M''}$  in Fus$^{simp}(\cV)$
and $(\cN^{gr}_{i})_{j,k}={\cN^{gr}_{M_i,M_j}}^{M_k}$ in Fus$^{gr}(\cV)$.
As before, they define representations of Fus$^{simp}(\cV)$ respectively Fus$^{gr}(\cV)$.

\begin{itemize}\item[\textbf{(V1)$'$}]  There are matrices $S^{simp\,\otimes}$ resp. $S^{gr\,\otimes}$
which simultaneously put the $\cN_\lambda$ resp. $\cN^{gr}_i$ into  block diagonal form: e.g.
$$S^{gr\,\otimes-1} \cN^{gr}_i  S^{gr\,\otimes}=\mathrm{diag}(B^{gr\,1}(M_i),B^{gr\,2}(M_i),\ldots,B^{gr\,r}(M_i))$$
where for each $h$, $M_i\mapsto B^h(M_i)$ defines
an indecomposable representations of the Grothendieck ring Fus$^{gr}(\cV)$ (and similarly for Fus$^{simp}(\cV)$).
\end{itemize}

(V1)$'$ is not deep: the fusion matrices $\cN^{gr}_i$ form a representation of Fus$^{gr}(\cV)$,
so $S^{gr\,\otimes}$ is the change-of-basis which decomposes that representation into a direct
sum of irreducibles. 
In \cite{R,PRR},  the  matrices $\cN^{simp}_{M}$ and  ${\cN^{gr}_{i}}$ 
for the $\cW_p$ models have been explicitly block-diagonalized: we review this work in Section 
3.3.2. We also should block-diagonalise the full tensor ring Fus$^{full}(\cV)$, in the sense of \cite{CR1, CR2}. To our knowledge this has never been explored.

To our knowledge, no analogue of (V2)$'$ has been explored. The Hopf link invariants of (V2) are still defined, but vanish on projective modules (at least if the category is rigid). 
If they did not vanish they would no longer yield tensor  or Grothendieck ring representations. The reason for this is that dim$\,\mathrm{End}(M_\lambda)$
can be larger than one. We propose a general  (V2)$'$ in Section 3.1.2.

Nothing general about 
{(V3)$'$} is known or has been explicitly conjectured. However,  Lyubashenko \cite{L} explained that the modular group acts also on certain 
nonsemi-simple braided tensor categories.  \cite{FGST1} observed that the $\text{SL}(2, \mathbb Z)$ action on the center of the restricted quantum group $\Uq$ of $sl_2$ at $2p$-th root of unity $q$,  coincides with the one on the space of trace and pseudo-trace functions of the $\cW_p$-triplet algebra.

The categorical SL$(2,\bbZ)$-action discussed last subsection extends here as follows. Write $\cC=\mathrm{Mod}^{g.r}(\cV)$.
According to \cite{Sh}, the
coend/Hopf algebra/Frobenius algebra/adjoint algebra $\cH$ is simply $UR(\cV)$, where $U:Z(\cC)\rightarrow\cC$ is
the forgetful functor from the Drinfeld centre of $\cC$, and $R$ is its right adjoint. The space of formal 1-point functions on the
torus is Hom$_\cV(\cV,\cH)$. The ordinary characters \eqref{char} span a subspace (of dimension equal to
the number of simple $\cV$-modules). We turn to this next.

\subsection{Miyamoto's modularity theorem}\label{sec:Miy}

Let $\cV$ be strongly-finite. Miyamoto \cite{M1} copies from Zhu \cite{Z} the definition of
\textit{formal 1-point functions}. These are, among other things, complex-valued functions $F(f,\tau)$ where
$\tau\in\bbH$ and $f=f(\tau)$ is a finite combination of vectors in $\cV$ with coefficients which are
modular forms of SL$(2,\bbZ)$ (so polynomials in $E_4(\tau)$ and $E_6(\tau)$).
Let $\cF(\cV)$ denote the space of all formal 1-point functions.

 It is easy to show directly from the definition that
$\cF(\cV)$ carries an action of SL$(2,\bbZ)$ through M\"obius transformations on $\tau\in\bbH$ as usual. 
The main result in \cite{Z} is that, for strongly-rational $\cV$, this space $\cF(\cV)$ has a 
basis over $\bbC[E_4,E_6]$ given by the characters $\chi_M(v,\tau)$ in \eqref{char}.
The
main result (Theorem 5.5) of \cite{M1} is that, when $\cV$ is strongly-finite,  $\cF(\cV)$ is finite-dimensional, and spanned by the pseudo-trace functions
$S^{M,\phi}(u,\tau)$, where $M$ is a `generalized Verma module interlocked with a symmetric linear functional $\phi$'
of the $n$th Zhu algebra $A_n(\cV)$ for $n$ sufficiently large (we'll try to explain this shortly).
These functions $S^{M,\phi}(u,\tau)$ are quite
difficult to compute in practice. Even the dimension of $\cF(\cV)$ is hard to determine in practice.
Although Miyamoto's approach is natural from the associative algebra point-of-view, it is not
at all so from the quantum field theory point-of-view. (We address this directly in Section \ref{sec:char} below.)
This makes it hard for many researchers to understand, despite its obvious importance. 
For this reason, in this subsection we'll supply some of the background and motivation for
Miyamoto's work.

 Let $S^{(j)}$ denote the basis of 
 $\cF(\cV)$ found by Miyamoto. Fix $u\in \cV$ with $L[0]u=ku$, then the vector $\bbX(\tau)$ with components $\bbX(\tau)_j=S^{(j)}(u,\tau)$ is a weight-$k$ vector-valued modular form.
Taking $u=\mathbf{1}$, the vacuum, recovers what could be called graded-pseudo-dimensions,
e.g. the familiar graded-dimensions $\chi_M(\tau)=q^{h_M-c/24}\sum_{n=0}^\infty\textrm{dim}\,M_nq^n$.

 Incidentally, the presence of $\bbC[E_4,E_6]$ in the arguments of formal 1-point functions is not significant.
  But the presence there of $u\in \cV$ is vital. It is common in the literature to drop $u$ from
  the character \eqref{char} and
  consider only the graded-dimensions $\chi_M(\tau)$. The reason for considering 1-point functions is  to avoid accidental
 linear dependencies between characters. For example, a module and its contragredient will always have the same graded-dimensions,
 but usually different characters (1-point functions).

Here is the problem in the strongly-finite case. Choose any indecomposable $\cV$-module $M$, and let $M({n})$ ($n=0,1,\ldots$) denote the generalized $L_0$-eigenspace of eigenvalue $h+n$ ($h$ is the conformal weight of $M$).
On each subspace $M({n})$,  $L_0-c/24$ acts as $(n+h-c/24)\mathrm{Id}+L_{nil}$ where $L_{nil}^k=0$ for some $k$ independent
of $n$. We obtain  for the usual character
\begin{equation}\chi_M(v,\tau)=\sum_{n=0}^\infty\mathrm{Tr}|_{M(n)}(o(v)q^{L_0-c/24})=\sum_{n=0}^\infty 
q^{n+h-c/24}\sum_{j=0}^{k-1}\frac{(2\pi\Ii\tau)^j}{j!}\mathrm{Tr}|_{M(n)}(o(v)\,L_{nil}^j)\,.\label{FMordinary}\end{equation}
It's not hard to show that on $M(n)$ the nilpotent operator $L_{nil}$ commutes with all operators $o(v)$. This means 
for any $j>0$ and any $v\in\cV$, $(o(v)\,L_{nil}^j)^k=o(v)^kL_{nil}^k=0$. But the trace of a nilpotent operator
is always zero,  so \eqref{FMordinary} collapses to
\[\chi_M(v,\tau)=\sum_{n=0}^\infty \mathrm{Tr}|_{M(n)}(o(v))q^{n+h-c/24}\,.\]
 Thus the ordinary trace can never see the nilpotent part of $L_0$, and $\chi_M(v,\tau)=\chi_N(v,\tau)$ whenever
$\cV$-modules  $M,N$ have the same composition factors. 

Nor can we obtain any more 1-point functions if we insert an endomorphism $f$ of $\cV$-module $M$: by definition,
$f$ will commute with all zero-modes and hence with $L_0$, so $f\circ o(v)\circ L_{nil}$ is still nilpotent. 

What we need is a way to generalize the trace
 of $M(n)$-endomorphisms so that the terms $o(v)L_{nil}^j$ with $j>0$ can contribute. We need to see
 the off-diagonal parts of the endomorphisms.

There is a classical situation where this happens. Let $A$ be a finite-dimensional  associative algebra with 1
over $\bbC$. A linear functional $\phi:A\rightarrow\bbC$ is called  \textit{symmetric}
if $\phi(ab)=\phi(ba)$ for all $a,b\in A$. Let $SLF(A)$ denote the space of all symmetric linear functionals on $A$
--- it can be naturally  identified with the dual space $(A/[A,A])^*$.
Now suppose $W$ is a finitely-generated projective $A$-module. Then an $A$-coordinate system of $W$
consists of finitely many $u_i\in W$ and the same number of $f_i\in \mathrm{Hom}_A(W,A)$ such that $w=
\sum_if_i(w)u_i$ for any $w\in W$. Given an $A$-coordinate system, we can associate any endomorphism
 $\alpha\in\mathrm{End}_A(W)$ with a matrix $[\alpha]$
whose $ij$th entry is $[\alpha]_{ij}=f_i(\alpha(u_j))\in A$. Fix any symmetric linear functional $\phi\in SLF(A)$ and 
$A$-coordinate system $\{u_i,f_i\}$, and define the \textit{pseudo-trace} Tr$^\phi_W:\mathrm{End}_A(W)\rightarrow\bbC$
by \[\mathrm{Tr}^\phi_W(\alpha)=\phi(\mathrm{Tr}([\alpha]))=\sum_i\phi\left(f_i(\alpha(u_i))\right)\,.\]
The pseudo-trace $\phi_W$ is independent of the choice of $A$-coordinate system, and lies in $SLF(\mathrm{End}_A(W))$.
It satisfies Tr$^\phi_W(\alpha\circ \beta)=\mathrm{Tr}^\phi_V(\beta\circ\alpha)$ for any $\alpha\in\mathrm{Hom}_A(V,W),\beta\in
\mathrm{Hom}_A(W,V)$.

To apply this generalized notion of trace to our VOA setting, we need to find a finite-dimensional associative algebra $A$
and $\cV$-modules $M$ for which the generalized eigenspaces $M_n$ are projective $A$-modules.  In \cite{M1},
 $A$ is related to the $n$th Zhu algebra $A_n(\cV)$ for $n$ sufficiently large, and $M$ are certain $\cV$-modules.
\cite{AN} makes a different choice: it is elementary that  a $\cV$-module $M$  (and in fact all of its generalized eigenspaces  $M_{n}$) is a
module over the (finite-dimensional associative) algebra End$_\cV(M)$, so choose
 some
subalgebra $A$ of End$_\cV(M)$ so that $M$ is projective over $A$. 

In both of these (related) cases, we have a finite-dimensional associative algebra $A$ and a symmetric linear functional
$\phi\in SLF(A)$, and a (logarithmic) $\cV$-module $M$ such that each generalised $L_0$-eigenspace $M_{(n)}$ is
a finite-dimensional projective $A$-module. We can then define the \textit{pseudo-trace function} 
\[\chi^{\phi}_M(v,\tau)=\sum_{n=0}^\infty q^{h_M+n-c/24}\sum_{j=0}^{k-1}\frac{(2\pi \Ii\tau)^j}{j!}\mathrm{Tr}^\phi_{M({n})}(o(v)\,L_{nil}^j)\] 
This  is a 1-point function in $\cF(\cV)$.
 In the Miyamoto case, these span $\cF(\cV)$;
in the simpler Arike--Nagatomo case, there is no such theorem.

Let's make this more concrete with an example arising in the $\cW_p$ models.
Suppose there are submodules $M,M'$ of a projective module $P$
such that 
$P/M'\cong M$, as $\cV$-modules. Suppose in addition that $M$ is a submodule of 
$M'$.   This happens for example with $M=0,M'=P$, but it also happens in $\cW_p$ with
 $M=\mathrm{soc}(P)$ and $M'=\mathrm{rad}(P)$ when $P=P^\pm_s$, $s\ne p$. Then in each generalized
eigenspace $P(n)$ of $L_0$, the operator $o(v)q^{L_0}$ has matrix form
\begin{equation}\label{MM'P}\left(\begin{matrix} A&B&C\\ 0&D&B^t\\ 0&0&A\end{matrix}\right)
\left(\begin{matrix} q^n&0&2\pi i q^n\tau\\ 0&q^n&0\\ 0&0&q^n\end{matrix}\right)\end{equation}
where the first row/column refers to the subspace $M$, the second row/column refers to any lift to $P$ of
$M'/M$, and the third row/column refers to any lift to $P$ of $P/M'$. Note that $C$ can be regarded
as an operator on $M$.
In this context, Miyamoto has a \textit{pseudo-trace} ptr$_{W(n)}$ of this operator, given by $2\pi i
\tau q^n\mathrm{Tr}\,C$. 
In particular, we see an off-diagonal part of the endomorphism.

More general 1-point functions, by inserting intertwiners, has been studied recently in this picture \cite{Fi}.

\subsection{Braided tensor categories}\label{sec:btc}

The books \cite{T}, \cite{EGNO}  are good references on tensor categories. Some features for nonsemi-simple ones relevant for us are included in \cite{CGP, GKP1}. 
Let $\mathcal C$ be our tensor category. We assume it to be strict so that we don't have to worry about associativity isomorphisms. This comes with a warning, the tensor category of a VOA is not strict, meaning that associativity isomorphisms are not trivial. One thus needs to always verify that theorems also hold in the non-strict case. A very useful theorem is here that every braided tensor category is braided equivalent to a free braided tensor category and in the latter morphisms composed out of braidings and associativity only depend on its braid diagram, see \cite{JS}.

The category is called \textit{rigid} if for each object $M$ in the category, there is a \textit{dual}
$M^*$ and morphisms $b_V\in\mathrm{Hom}(1,V\otimes V^*)$ (the co-evaluation) and
$d_V\in\mathrm{Hom}(V^*\otimes V,1)$ (the evaluation) such that
\begin{equation}\nonumber
 (\mathrm{Id}_V\otimes d_V)\circ(b_V\otimes \mathrm{Id}_V)=\mathrm{Id}_V\,,\ \
(d_V\otimes\mathrm{Id}_{V^*})\circ(\mathrm{Id}_{V^*}\otimes b_V)=\mathrm{Id}_{V^*}\label{rigid}
\end{equation}
Rigidity is required for the categorical trace and dimension, as we'll see shortly.
In the vertex operator algebra language, a morphism is an intertwiner and an object is a module. 
Given any modules $V,W$ in the category, the \textit{braiding} $c_{V,W}$ is an intertwiner in Hom$(V\otimes W,W\otimes V)$ satisfying 
\begin{equation}\label{braiding}
c_{U,V\otimes W}=(\mathrm{Id}_V\otimes c_{U,W})\circ(c_{U,V}\otimes \mathrm{Id}_W)\,,\ \
c_{U\otimes V,W}=(c_{U,W}\otimes \mathrm{Id}_V)\circ(\mathrm{Id}_U\otimes c_{V,W})
\end{equation}
and also for any intertwiners $f\in\mathrm{Hom}(V,V')$, $g\in\mathrm{Hom}(W,
W')$,
\begin{equation}\label{eq:braidnat}(g\otimes f)\circ c_{V,W}=c_{V',W'}\circ (f\otimes g)\,.\end{equation}
This property is called \textit{naturality} of braiding. 
It implies $c_{V,1}=c_{1,V}=1$ as well as the Yang-Baxter equation
$$(\mathrm{Id}_W\otimes c_{U,V})\circ(c_{U,W}\otimes\mathrm{Id}_V)\circ(\mathrm{Id}_U\otimes c_{V,W})
=(c_{V,W}\otimes\mathrm{Id}_U)\circ(\mathrm{Id}_V\otimes c_{U,W})\circ(c_{U,V}\otimes\mathrm{Id}_W)\,.$$
If the category $\cC$ is in addition  \textit{additive} (i.e. it has direct sums), then $(U\oplus V)\otimes W$ resp.\
$W\otimes(U\oplus V)$ are isomorphic with $U\otimes W\oplus V\otimes W$ resp.\ $W\otimes U\oplus W\otimes V$,
and using these isomorphisms the naturality \eqref{eq:braidnat} of the braiding
implies that we can make the identifications
\begin{equation}\label{eq:braidds}
c_{U\oplus V,W}=c_{U,W}\oplus c_{V,W}\,,\ \ c_{W,U\oplus V}=c_{W,U}\oplus c_{W,V}\,.\end{equation}
In a (semi-simple) modular tensor category, the braiding is needed for the topological $S$-matrix, as we'll see
shortly.

Given any module $V$ in our category, the \textit{twist} $\theta_V\in\mathrm{Hom}(V,V)$ satisfies,
for any module $V$ and intertwiner $f\in\mathrm{Hom}(V,V)$, 
\begin{equation}\nonumber
\theta_{V\otimes W}=c_{W,V}\circ c_{V,W}\circ (\theta_V\otimes \theta_W)\,,\ \
\theta_V\circ f=f\circ\theta_V\label{twist}
\end{equation}
This implies $\theta_1=1$. The twist directly gives us the (diagonal) topological $T$-matrix in a modular tensor 
category: for any simple object $V$, $\theta_V$ is a number and $T_{V,V}= \theta_V$.
Define morphisms
\[
b'_V:=\left(\mathrm{Id}_{V^*}\otimes \theta_V\right)\circ c_{V, V^*}\circ b_V \ \in \ \mathrm{Hom}\left(1, V^*\otimes V\right)
\]
and
\[
d'_V:= d_V\circ c_{V, V^*}\circ\left(\theta_V\otimes \mathrm{Id}_{V^*}\right) \ \in \ \mathrm{Hom}\left(V\otimes V^*, 1\right).
\]
Given any intertwiner $f\in\mathrm{Hom}(V,V)$, the \textit{(categorical- or quantum-)trace} is defined by 
\[
\text{tr}(f)=d'_V\circ\left(f\otimes \mathrm{Id}_{V^*}\right)\circ b_V\ \in\ \mathrm{ Hom}(1,1)\,.
\]
We can identify $\mathrm{ Hom}(1,1)$ with $\bbC$. Diagramatically, an endomorphism $f$
is a vertical arrow and its trace corresponds to attaching top and bottom into a circle.
The trace satisfies tr$(f\circ g)=\mathrm{tr}(g\circ f)$ and tr$(f\otimes g)=\mathrm{tr}(f)\mathrm{tr}(g)$.
We define the \textit{(categorical- or quantum-)dimension} to be dim$(V)=\mathrm{tr}(\mathrm{Id}_V)$. It will be positive
if the category is unitary.
Computations are greatly simplified using a graphical calculus, some basic notation is 
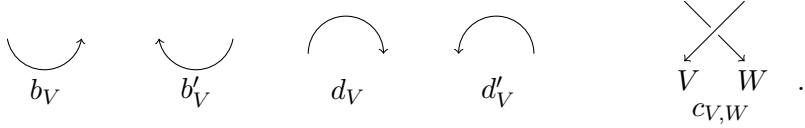
\begin{center}
\begin{tikzpicture}

\draw[br] (0, 3) arc (190:350:.5cm);
\node[text width=3cm] at (1.8, 2.3)    {$b_V$};

\draw[br] (3, 3) arc (350:190:.5cm);
\node[text width=3cm] at (3.8, 2.3)    {$b'_V$};

\draw[ar] (5, 2.8) arc (0:180:.5cm);
\node[text width=3cm] at (5.8, 2.3)    {$d_V$};

\draw[br] (7, 2.8) arc (0:180:.5cm);
\node[text width=3cm] at (7.8, 2.3)    {$d'_V$};

\draw [ar](9,2.7) -- (9.8,3.5);
\draw (9,3.5) -- (9.35,3.15);
\draw [br](9.45,3.05) -- (9.8,2.7);

\node[text width=3cm] at (10.4, 2.45)    {$V$};
\node[text width=3cm] at (11.2, 2.45)    {$W\ \ \ .$};
\node[text width=3cm] at (10.6, 2.0)    {$c_{V,W}$};

\end{tikzpicture}\end{center}
Details on the graphical calculus are found in the textbooks \cite{T, K, BK}.

By a \textit{finite tensor category} \cite{EO} we mean an abelian rigid tensor category over $\bbC$, where morphism
spaces are finite-dimensional, there are only finitely many simple objects, each of which has  a projective
cover, every object has finite length in the sense of  Jordan--H\"older,  and  the endomorphism algebra
of the tensor unit is $\bbC$.
By a \textit{fusion category} (see e.g. \cite{EGNO}) we mean a finite tensor category which is in addition semi-simple.
By a \textit{ribbon category} \cite{T} we mean a strict braided tensor category equipped with left duality ${}^\vee X$ and a twist $\theta_X\in\mathrm{End}(X)$.

In a ribbon finite tensor category, the categorical $S$-matrix is defined by $$S^\hopflink_{U,V}=
\mathrm{tr}(c_{U,V}\circ c_{V,U})\in \bbC,$$ for any indecomposable modules $U,V$. 
Its graphical representation is the Hopf link \eqref{hopflink}.
By a \textit{modular tensor category} \cite{T} we mean a fusion category which is ribbon, with invertible matrix $S^\hopflink$.

The categorical $S$-
and $T$-matrices define a projective representation of SL$(2,\bbZ)$. It satisfies $S^\hopflink_{U,V}
=S^\hopflink_{V,U}$, $S^\hopflink_{1,V}=\mathrm{dim}(V)$, as well as \textit{Verlinde's formula}
\begin{equation}
\nonumber
\sum_W \cN_{U,V}^WS^\hopflink_{W,X}=(\mathrm{dim}\,X)^{-1}S^\hopflink_{U,X}S^\hopflink_{V,X}\label{verlinde}
\end{equation}
Ideally (e.g. when the VOA is strongly-rational), the categorical $S$- and $T$-matrices will
agree up to scalar factors with the modular $S$- and $T$-matrices defined through the
VOA characters.

Return for now to a finite tensor category $\cC$. Assume it is spherical.
By a \textit{negligible morphism} $f\in\mathrm{Hom}_\cC(U,V)$ we mean one for which the categorical trace Tr$_\cC(g\circ f)=0$ 
for all $g\in\mathrm{Hom}_\cC(V,U)$. The negligible morphisms form a subspace $I_\cC^{neg}(U,V)$ of Hom$_\cC(U,V)$, closed under taking duals, arbitrary compositions,
 as well as  arbitrary tensor products. Consider the category $\overline{\cC}$ whose objects are the same as
 those of $\cC$, but whose Hom-spaces are Hom$_\cC(U,V)/I_\cC^{neg}(U,V)$. This category was originally defined in
 \cite{BW}; see also p.236 of \cite{EGNO}. If $U,V$ are indecomposable in $\cC$ and $f\in \mathrm{Hom}_\cC(U,V)$ is not 
 an isomorphism, then $f$ is negligible. Then $\overline{\cC}$ is a spherical tensor category, whose simple objects
 are precisely the indecomposable objects of $\cC$ with nonzero categorical dimension (the indecomposables with
 dimension 0 are in $\overline{\cC}$ isomorphic to the 0-object). Moreover, two simple objects in $\overline{\cC}$ are
 isomorphic iff they are isomorphic indecomposables in $\cC$. So $\overline{\cC}$ will generally have infinitely
 many inequivalent simples, but it is semi-simple in the sense that every object in $\overline{\cC}$ is a direct sum
 of simples. By the \textit{semi-simplification} $\cC^{ss}$ of $\cC$, we mean the full subcategory of $\overline{\cC}$
 generated by the simples of $\cC$. Then $\cC^{ss}$ will also be a semi-simple spherical tensor category. If in $\cC$ the
 tensor product of simples is always a direct sum of simples and projectives (as it is for $\cW_p$ and the symplectic
 fermions), then $\cC^{ss}$ will be a fusion category, whose Grothendieck ring is the quotient of that of $\cC$ by
 the modules in $\cC$ which are both simple and projective. If $\cC$ is ribbon, so is both $\overline{\cC}$ and
 $\cC^{ss}$. 

\begin{remark} 

Note though that the semi-simplification of a nonsemi-simple finite tensor category $\cC$ can never be unitary,
as some quantum-dimensions must be negative in order that the projectives of $\cC$ have 0 quantum-dimension.

\end{remark}

\subsection{Lyubashenko's modularity and the coend}\label{lyb}

Lyubashenko realized how modular group and more general mapping class group representations also arise in the context of certain non semisimple ribbon categories. 
The important object here is the coend of $\cC$. 
It is denoted by
\[
\coend = \int^{X\in\;\cC} X\otimes X^\vee
\]
and defined by \cite{L3}
\[
\bigoplus_{f : A \rightarrow  B \in \cC} A\otimes B^\vee   \xrightarrow{\ f\otimes B^\vee - A \otimes f^t  \ }  \bigoplus_{X\in\;\cC} X\otimes X^\vee \longrightarrow \coend \longrightarrow 0 
\]
where $f^t: B \rightarrow A^\vee$ is the transposed of $f:A\rightarrow B^\vee$.
Lyubashenko defines $2$-modular categories \cite{L} (see also Definition 1.3.1 of \cite{L3}). These are abelian ribbon categories with additional axioms.  
They imply the existence of a morphism $\mu : \one \rightarrow \coend$ (called an integral) solving an equation relating open Hopf links dressed with twist to the inverse of the twist. This is called the quantum Fourier transform and explicitely stated in Section 1 of \cite{L3}. The existence of the integral $\mu$ then allows to show that a morphism $S : \coend \rightarrow \coend$ defined via the integral and the monodromy together with another morphism $T$ that is defined via the twist define a projective representation of SL$(2;\mathbb Z)$. More generally there is even a projective representation of the mapping class group of genus $g$ with $n$-punctures on Hom$\left(X_1 \otimes \dots \otimes X_n, \coend^{\otimes g}\right)$. Clearly one expects and would like to understand a relation of these to conformal blocks. 

The role of the coend in conformal field theory has been investigated very recently in \cite{FS3}. Very nice and probably also useful properties are that the coend has the structure of a Hopf  algebra (more precisely a Hopf monad) and there is a braided equivalence between the category of $\coend$-modules in $\cC$ and the center $Z(\cC)$ of $\cC$ \cite{DS, BrV}.

\subsection{Example 1: The triplet algebra $\cW_p$}

Our main reference on the triplet algebra is \cite{TW}.
The triplet models  $\cW_p$  ($p\ge 2$ an integer) are a family of VOAs with central charge $c_p=1-6(p-1)/p$.
They are logarithmic and strongly-finite \cite{AM}. As mentioned in Section 2.4, this implies they have finitely many
simple (ordinary) modules, whose (pseudo-)characters  form a vector-valued modular function
for SL$(2,\bbZ)$, and which yields a  braided tensor category Mod$^{g.r}(\cW_p)$. 
$\cW_p$ has been studied by Adamovic and Milas in a series of papers \cite{AM, AM3, AM4, Mil2}; in \cite{AM4} Zhu's algebra was described and in \cite{AM3, Mil2} it was shown that the space of 1-point functions is $3p-1$-dimensional.
The tensor product and rigidity for $\cW_p$ were determined
in \cite{TW} (the tensor product had been previously conjectured in \cite{FHST},\cite{GR}).

The simple modules are $X^+_s$ and $X^-_s$ for $1\le s\le p$. The tensor unit (vacuum) is $X^+_1$.
Their
projective covers are $P^+_s$ and $P^-_s$, where $P^\pm_p=X^\pm_p$ and 
$$0\rightarrow Y_s^+\rightarrow P_s^+\rightarrow X^+_s\rightarrow 0,\qquad
0\rightarrow Y_s^-\rightarrow P_{p-s}^-\rightarrow X^-_{p-s}\rightarrow 0$$
for $1\le s\le p-1$, where $Y^\pm_s$ denotes the reducible but indecomposable modules
$$0\rightarrow X^+_s\rightarrow Y^+_s\rightarrow  2{\cdot}X^-_{p-s}\rightarrow 0\,, \qquad
0\rightarrow X^-_{p-s}\rightarrow Y^-_s\rightarrow  2{\cdot} X^+_{s}\rightarrow 0\,.$$
These $4p-2$ irreducible and/or projective $\cW_p$-modules $X^\pm_s,P^\pm_s$  are the most important ones. They are closed
under tensor product and form the subring Fus$^{simp}(\cW_p)$.

\cite{NT} proved the equivalence  as abelian categories of Mod$^{g.r}(\cW_p)$ with the category Mod$^{fin}(\Uq)$ of
finite-dimensional modules of the restricted quantum group  $\Uq$ at $q=e^{\pi\Ii/p}$, and both \cite{FGST2, KS} determined all
indecomposables in Mod$^{fin}(\Uq)$ (although the key lemma, describing pairs of matrices up to simultaneous 
conjugation, is really due to Frobenius (1890)). Hence:
\begin{theorem}\label{thm:Wpindec} 
(\cite{TW})  The complete list of indecomposable grading-restricted weak $\cW_p$-modules, up to equivalence, is:
\begin{enumerate}
\item  for each $1\le j\le p$ and each sign, the simple modules $X^\pm_j$ and their projective covers  $P_j^\pm$
($P_p^\pm=X^\pm_p$);
\item  for each $1\le j< p$, each sign, and each $d\ge 1$,  $G^\pm_{j,d}$, $H^\pm_{j,d}$, and $I^\pm_{j,d}(\lambda)$, where  $\lambda\in\bbC\bbP^1$.
\end{enumerate}
\end{theorem}

Consider the   full subcategory relevant to our discussion, whose  objects are finite direct sums of $X^\pm_s,P^\pm_s$ for $1\le s\le p$. 
 The fusion rules are completely determined from the following: 
\begin{equation}\label{eq:minsetfus}\nonumber
\begin{split}
X^{\epsilon}_1\otimes X^{\epsilon'}_s&=X^{\epsilon\epsilon'}_s\,,\ \ \ X_1^\epsilon\otimes P^{\epsilon'}_s=P_s^{\epsilon\epsilon'}\,,\\
X^+_2\otimes X^+_s&=\left\{\begin{matrix}X^+_2&s=1\\ X^+_{s-1}\oplus X^+_{s+1}&2\le s<p\\
P^+_{p-1}&s=p\end{matrix}\right.\,,\\
X^+_2\otimes P^+_s&=\left\{\begin{matrix}P_2^+\oplus 2\cdot X^-_p&s=1\\ P^+_{s-1}\oplus P^+_{s+1}&2\le s<p-1\\
P^+_{p-2}\oplus 2\cdot X^+_p&s=p-1\end{matrix}\right. \,,\\
\end{split}
\end{equation}
together with associativity and commutativity. For example one can show
\[P^\epsilon_s\otimes P^{\epsilon'}_t=2{\cdot} X^+_s\otimes P_t^+\oplus 2{\cdot} X^-_{p-s}\otimes P_t^+\,.\]
In $\cW_2$, $X^+_2\otimes P^+_s$ must be replaced with $X^+_2\otimes P^+_1=2{\cdot}X^+_2\oplus 2{\cdot}X^-_2$.

As always, the $\bbZ$-span of the projective modules $\{P_s^{\pm}\}_{1\le s\le p}\cup\{X^\pm_p\}$ forms
an ideal in Fus$^{full}(\cW_p)$. The quotient of the tensor subring $\bbZ$-span$\{X^\pm_s,P^\pm_t\}$ by
that ideal is easily seen to be two copies of the fusion ring of $L_{p-2}\left(\mathfrak{sl}_2\right)$, more precisely it is isomorphic to the fusion ring of the rational VOA $L_{1}\left(\mathfrak{sl}_2\right)\otimes L_{p-2}\left(\mathfrak{sl}_2\right)$, where $X^-_1$ corresponds to the integrable highest weight module
$L(\Lambda_1)$ of $L_{1}\left(\mathfrak{sl}_2\right)$ and $X^+_s$ corresponds likewise to the module
$L((p-1-s)\Lambda_0+(s-1)\Lambda_1)$ of $L_{p-2}\left(\mathfrak{sl}_2\right)$ for $1\le s<p$. The proof is elementary,
obtained by comparing the tensor products of generators. In fact, more importantly, this also persists categorically:
the semi-simplification $(\mathrm{Mod}^{g.r}(\cW_p))^{ss}$ is a modular tensor category, namely the twist of that of
 $L_{p-2}\left(\mathfrak{sl}_2\right)
\otimes L_{1}\left(\mathfrak{sl}_2\right)$ by a simple-current of order 2 as discussed in \cite{CG}.

The characters (really graded-dimensions) of $X^\pm_s$ are
$$\mathrm{ch}[X^+_s](\tau)=\frac{1}{\eta(\tau)}\left(\frac{s}{p}\theta_{p-s,p}(\tau)+2\theta'_{p-s,p}(\tau)\right),\quad
\mathrm{ch}[X^-_s](\tau)=\frac{1}{\eta(\tau)}\left(\frac{s}{p}\theta_{s,p}(\tau)-2\theta'_{s,p}(\tau)\right)$$
where 
$$\theta_{s,p}(\tau,z)=\sum_{j\in \frac{s}{\sqrt{2p}}+\sqrt{2p}\bbZ}e^{2\pi \Ii\tau j^2/2}
e^{2\pi \Ii zj}$$
is the theta series associated to the coset $\frac{s}{\sqrt{2p}}+L\in L^*/L$ of the
even lattice $L=\sqrt{2p}\bbZ$ and 
$$
\theta_{s,p}(\tau)=\theta_{s,p}(\tau,0)\,,\qquad
\theta'_{s,p}(\tau)=\frac{1}{2\pi \Ii}\frac{\partial}{\sqrt{2p}\partial z}\theta_{s,p}(\tau,z)\Big\vert_{z=0}\,.
$$ 
The characters of $P^+_s$ and $P^-_{p-s}$ are both 
$$\mathrm{ch}[P^+_s](\tau)=\mathrm{ch}[P^-_{p-s}](\tau)=2\mathrm{ch}[X^+_s](\tau)+2\mathrm{ch}[X^-_{p-s}](\tau).$$

We define the pseudo-character to be the weight one part of characters (of course the choice of normalization is 
a convention)
\[
\mathrm{pch} [X^+_s] :=-4\Ii\tau \frac{\theta'_{p-s,p}\left(\tau\right)}{\eta\left(\tau\right)}\,,\qquad
\mathrm{pch} [X^-_s] := 4\Ii\tau \frac{\theta'_{s,p}\left(\tau\right)}{\eta\left(\tau\right)}\,.
\]
Not all characters are linearly independent and a basis is given by 
\[
\mathcal B:=\left\{ \ \mathrm{ch}[P^+_{\ell}], \mathrm{ch}[X^+_{\ell}], \mathrm{pch} [X^+_\ell], \mathrm{ch}[X^\pm_{p}]\ \Big| \ 1\leq \ell\leq p-1\ \right\}.
\]
Set $q:= e^{\frac{2\pi \Ii }{2p}}$.
The modular $S$-transformations is
\begin{equation}\nonumber
\begin{split}
\mathrm{pch} [X^+_s]\left(-\frac{1}{\tau}\right)&= 
\frac{4}{\sqrt{2p}}\sum_{\ell=1}^{p-1} (-1)^{\ell+s+p}\left(q^{\ell s}-q^{-\ell s}\right) \left( \mathrm{ch}[X^+_{\ell}]-\frac{\ell}{2p}\mathrm{ch}[P^+_{\ell}]\right)\\
\mathrm{ch} [X^+_s]\left(-\frac{1}{\tau}\right)
&=\frac{1}{\sqrt{2p}} \sum_{\ell=1}^{p-1} (-1)^{\ell+s+p+1} \left(q^{\ell s}-q^{-\ell s}\right) \mathrm{pch} [X^+_\ell]\\
\mathrm{ch} [P^+_s]\left(-\frac{1}{\tau}\right)
&=\frac{(-1)^{p}}{\sqrt{2p}}\Bigl((-1)^{p}\mathrm{ch} [X^+_p]+(-1)^{s}\mathrm{ch} [X^-_p]+ \sum_{\ell=1}^{p-1}(-1)^{\ell+s}\left(q^{\ell s}+q^{-\ell s}\right) \mathrm{ch}[P^+_{\ell}]\Bigr).
\end{split}
\end{equation}
We define the modular $S$-matrix coefficients via
\begin{equation}\nonumber
\begin{split}
\mathrm{ch}[V]\left(-\frac{1}{\tau}\right)&= \sum_{\text{characters}\ \text{in} \ \mathcal B} S^\chi_{V, W}\  \mathrm{ch}[W]   + \sum_{\text{pseudotraces}\ \text{in} \ \mathcal B} S^{\chi}_{V, W^0}\  \mathrm{pch}[W],\\
\mathrm{pch}[V]\left(-\frac{1}{\tau}\right)&= \sum_{\text{characters}\ \text{in} \ \mathcal B} S^\chi_{V^0, W}\  \mathrm{ch}[W]   + \sum_{\text{pseudochars}\ \text{in} \ \mathcal B} S^{\chi}_{V^0, W^0}\  \mathrm{pch}[W].
\end{split}
\end{equation}

\subsection{Example 2: Symplectic fermions}\label{SF}

Another class of strongly-finite VOAs is $SF^+_d$ for $d\ge 1$, the even part of the symplectic fermions.
The VOAs $SF^+_1$ and $\cW_2$ coincide, and actually the structure of $SF_d^+$ is inherited from $\mathcal W_p$ as $SF^+_d$ is an extension by an abelian intertwining algebra of the $d$-fold tensor product of $\mathcal W_2$.

$SF_1=SF$ is a vertex operator super algebra. It has just itself as its only simple module, and the projective cover $P$ of it has the form 
$0\rightarrow T\rightarrow P \rightarrow T\rightarrow 0$, where $T$ satisfies $0\rightarrow SF \rightarrow T\rightarrow SF \rightarrow 0$ (both non-split). 
$SF_d$ is its $d$-fold tensor product, it thus has also only one simple module and the projective cover is the tensor product of the projective covers of the components. 

$SF^+_d$ has exactly 4 simple modules \cite{Abe}: $SF^\pm_d$ and $SF^\pm(\theta)_d$, where $SF^-_d$ is the odd
part of the symplectic fermion vertex operator superalgebra, and  $SF^\pm(\theta)_d$ are the even/odd parts of a
twisted module of the symplectic fermions.
$SF_1^+$ is just $\mathcal W_2$ and its fusion is thus (for $\epsilon, \nu\in\{\pm\}$)
\begin{equation}\nonumber
\begin{split}
SF^\epsilon_1 \otimes SF^\nu_1 &= SF_1^{\epsilon\nu},\qquad\qquad\qquad\qquad
SF^\epsilon_1 \otimes SF^\nu_1(\theta) = SF_1^{\epsilon\nu}(\theta)\,, \\
SF^\epsilon_1(\theta) \otimes SF^\nu_1 (\theta)&= P_1^{\epsilon\nu},\qquad\qquad\qquad\qquad\qquad \ \,\,
SF^\epsilon_1 \otimes P^\nu_1 = P_1^{\epsilon\nu}\,, \\
SF^\epsilon_1(\theta) \otimes P^\nu_1 &= 2{\cdot}SF^+_1(\theta)\oplus 2{\cdot}SF^-_1(\theta),\qquad\ \
P^\epsilon_1 \otimes P^\nu_1 =  2{\cdot}P^+_1\oplus 2{\cdot}P^-_1\,. 
\end{split}
\end{equation}
Here $P^\pm_1$ are the even/odd part of $P$.
The modules $SF^\pm_d, SF^\pm(\theta)_d$ and $P^\pm_d$ are
the even/odd part of the $d$-fold tensor product of the corresponding $SF$-modules, that is they
are the following $\left(SF_1^+\right)^{\otimes d}$-modules:
\begin{equation}\nonumber
\begin{split}
X^+_d &=\bigoplus_{\substack{\epsilon\in \mathbb F_2^d\\ \epsilon^2=0 }} \bigotimes_{i=1}^d X^{\epsilon_i}_1, \qquad
X^-_d =\bigoplus_{\substack{\epsilon\in \mathbb F_2^d\\ \epsilon^2=1 }} \bigotimes_{i=1}^d X^{\epsilon_i}_1,  \qquad X\in \{ SF, SF(\theta), P\}.
\end{split}
\end{equation}
Their tensor ring can be deduced from the one of $\left(SF_1^+\right)^{\otimes d}$ and is isomorphic to the case of $d=1$.
For details on this procedure see section 5 of \cite{AA}. However, note that the tensor ring is not the Klein four group. 
The symplectic fermions $SF^+_d$ are rigid (Proposition 3.23 of \cite{DR}).
 The characters (really, graded-dimensions) of the simple modules are
 \begin{equation}\nonumber
\begin{split}
\mathrm{ch}[SF^\pm_d](\tau)&=\frac{1}{2}\left(\frac{\eta(2\tau)^{2d}}{\eta(\tau)^{2d}}\pm\eta(\tau)^{2d}\right),\\ 
\mathrm{ch}[SF^\pm(\theta)_d](\tau)&=\frac{1}{2}\left(\frac{\eta(\tau)^{4d}}{\eta(2\tau)^{2d}\eta(\tau/2)^{2d}}\pm \frac{\eta(\tau/2)^{2d}}{\eta(\tau)^{2d}}\right).
\end{split}
\end{equation}
The conformal weights are $0,1,-d/8,(4-d)/8$ respectively. We see that $SF^\pm(\theta)_d$ are projective in addition to
being simple.
Abbreviate these characters as $\chi^\pm,\chi^\pm_\theta$.
Modularity of the characters of the simple modules is given by:
\begin{equation}\nonumber
\begin{split}
\chi^\pm(\tau+1)&=e^{\pi \Ii d/6}\chi^\pm(\tau)\,,\\
 \chi^\pm(-1/\tau)&=\pm\frac{(-\Ii\tau)^d}{2}\left(\chi^+(\tau)
-\chi^-(\tau)\right)+\frac{1}{2^{d+1}}
\left(\chi_\theta^+(\tau)-\chi_\theta^-(\tau)\right)\,,\\
\chi^\pm_\theta(\tau+1)&=\pm e^{-\pi \Ii d/12}\chi^\pm(\tau)\,,\\ 
 \chi_\theta^\pm(-1/\tau)&=\pm 2^{d-1}\left(\chi^+(\tau)+\chi^-(\tau)\right)+\frac{1}{2}
\left(\chi_\theta^+(\tau)+\chi_\theta^-(\tau)\right)\,.\\ \end{split} \end{equation}
Two indecomposable but nonsimple $SF^+_d$-modules $\widehat{SF}_d^\pm$  are constructed in \cite{Abe}. They both have characters $2^{2d-1}\eta(2\tau)^{2d}/\eta(\tau)^{2d}$ and socles $SF^\pm_d$. They have Jordan--H\"older length $2d$ with composition factors
$ SF^\pm_d, SF^\mp_d\otimes \bbC^{2d},\ldots, SF^{(-1)^r}_d\otimes \bbC^{({2d\atop r})},\ldots, SF^\pm_d$ \cite{Abe}.
Thus their images in the Grothendieck ring are both $2^{2d-1}[SF^+_d]+2^{2d-1}[SF^-_d]$.
$L_0$ for both these indecomposable modules acts by Jordan blocks of size $d+1$ \cite{AN}. 
\cite{AN} show that the dimension of the SL$(2,\bbZ)$-representation is at least $2^{2d-1}+3$, and conjecture
this is the exact dimension (since $SF^+_1$ coincides with $\cW_2$, we know this conjecture is correct for $d=1$). 

$\widehat{SF}^\pm_d$ are projective $SF_d$-modules and hence the projective covers of $SF^\pm_d$. This follows for example by viewing $SF^+_d$ as an algebra for $\left(SF^+\right)^{\otimes d}$ according to \cite{HKL}, then $\widehat{SF}^\pm_d$ is the induced module of a projective $\left(SF^+\right)^{\otimes d}$-module and as such projective itself as an $SF^+_d$-module, see \cite{EGNO, CKM}. We will explain this in detail in \cite{CG}.

Symplectic fermions, having automorphism group the symplectic Lie group,  allow for many other orbifolds \cite{CL3} than 
just by $\mathbb Z_2$  and any orbifold of a $C_2$-cofinite VOA by any solvable finite group is itself $C_2$-cofinite 
\cite{Miy4}. Thus giving rise to many more examples of strongly-finite VOAs.

\section{Our results}\label{sec:results}

\subsection{Conjectures, observations and questions}\label{sec:conj}

Let $\cV$ be a strongly-finite VOA, i.e. a simple $C_2$-cofinite  VOA of CFT-type with $\cV$ isomorphic to $\cV^\vee$.
Based on the preceding remarks, one could hope that the following are true, at least in broad strokes,
though surely some details will have to be modified.
It is perhaps premature to make precise conjectures, but we propose the following as a way to help
guide future work.

The whole area needs more independent examples, to probe conjectures such as those collected
below, and also so suggest new conjectures. We address this challenge in Section 3.4 below.

\subsubsection{Characters and pseudo-characters}\label{sec:char}

A problem is that, although the generalised trace functions of Miyamoto are very natural from the point of view of
associative algebras, they are not natural from the point of view of VOAs. In quantum field theory, one should be
able to recover all $n$-point functions from ordinary traces with insertions. In particular, we should be able
to recover all modular forms in Miyamoto's space $\cF(\cV)$, through the 1-point functions $\cZ^\cY(v,\tau)$ of \eqref{1ptfn},
where $\cY$ is some intertwiner. Let's make this idea more explicit.

 Write $P^{1}$ for the unique
 projective $\cV$-module whose socle is $\cV$, i.e. $P^1$ is indecomposable and contains $\cV$
 as a submodule. 
By Conjecture \ref{unimodular} below, $P^1$ should be the projective cover of $\cV$.
There is only one submodule of $P^1$ isomorphic to $\cV$, but typically $P^1$ will have several
other composition factors (i.e. subquotients) isomorphic to $\cV$; by a \textit{$\cV$-lift} we mean
a subspace $V$ of some submodule $M$ in $P^1$ such that $M/N$ is isomorphic to $\cV$ 
for some submodule $N$ of $M$, and $V\cap N=0$. 

\begin{conjecture}  {Let $\tilde{\cF}(\cV)$ be the span over $\bbC[E_4(\tau),E_6(\tau)]$ of all 
functions $\tau\mapsto Z^\cY(u,\tau)$, as $\cY$ runs over any intertwiner of type $\left({P\atop P^1,P}
\right)$ for $P$ any projective cover of a simple $\cV$-module, and any $u$ in any $\cV$-lift.
Then  $\tilde{\cF}(\cV)$ equals the space of all $F\in\cF(\cV)$ regarded as a function of $\tau$ only.}
\end{conjecture}

The space $\cF(\cV)$ is finite-dimensional over $\bbC[E_4(\tau),E_6(\tau)]$ when the functions $F$ are regarded as functions of
2 variables $v,\tau$, but $\tilde{\cF}(\cV)$ is infinite-dimensional, so information is lost.
We would like to strengthen the conjecture: e.g. show that there are lifts to $P^1$ of each subquotient of $P^1$ isomorphic to $\cV$, such that the 1-point functions restricted to those lifts are formal 1-point
functions in  Miyamoto's space $\cF(\cV)$.

Ordinary intertwining operators $\cY(v,z)$ (i.e. without any log$\,z$ terms) will commute with the nilpotent operator $L_0^{nil}=L_0-L_0^{ss}$, and so their 1-point functions $Z^\cY_P(v,\tau)$
will be pure $q$-series by an argument given in Section \ref{sec:Miy}. But logarithmic intertwiners can
have log$\,q$ contributions, as shown explicitly by Kausch \cite{Kau} some time ago (see also
Section 7.1 of \cite{Ru}). 

This conjecture reduces to Zhu's theorem when $\cV$ is strongly-rational.
A modular group action on the 1-point functions $Z^\cY(v,\tau)$ in the strongly-finite context
is studied in \cite{Fi}. If one is interested only in weight-0 modular forms (generalised graded-dimensions) then one should restrict to $v$ of minimal generalised conformal weight.
We can expect, from \cite{M0} (which was proved in the strongly-rational
context), that  choosing $v$ with $L_{[0]}v=kv$ will get modular forms of weight $k$.

Intertwiners are intimately related to tensor products. Suppose $P$ is the 
projective cover of some simple module $M_k$. Then its tensor product with $P^1$  will be
$$P^1\otimes P\cong \sum_j \langle P^1,M_j\rangle\, M_j\otimes M_k$$
where we sum over all simple modules $M_j$, and where $\langle P^1,M_j\rangle$ denotes
the multiplicity of $M_j$ as a composition factor of $P^1$. In other words, we get an independent
intertwiner $\cY$ of type $\left({P\atop P^1,P}\right)$, every time $\cV$ is a composition factor of $P^1$. The intertwiner associated to a subquotient $M/N\cong\cV$ in $P^1$ will vanish in $N$.

Miyamoto's space $\cF(\cV)$ for the symplectic fermions is surprisingly large: of dimension at least
(and probably exactly) $2^{2d-1}+3$. But each of the 4 possible indecomposable $P$ will
have $2^{2d-1}$ intertwiners (since $P^1$, the projective cover of $\cV$, has precisely
$2^{2d-1}$ composition factors of type $\cV$). 

Let's make this more concrete by considering the $\cW_2$ model. This is of the type discussed in Section 2.4 (recall \eqref{MM'P}). Take $P=P^1=P^+_1$. Then $M=X^+_1$ and $M'$ has socle
$X_1^+$ and top $X^-_1\oplus X^-_1$. 
There are 2 linearly independent intertwiners of type $\left({P\atop P,P}\right)$, corresponding to the
2 $X^+_1$'s in $P^1$: the submodule soc$(P^1)$ and the quotient top$(P^1)$. Choosing $v$ in the socle for the first intertwiner corresponds to $o_\cY(v)$ a nonzero scalar multiple of $o_M(v)$, so 
to  an endomorphism $o_\cY(v)q^{L_0}$ given in \eqref{MM'P}; the corresponding trace
$Z_\cY(v,\tau)$ will be a constant multiple of the ordinary 1-point function $2\mathrm{ch}[{X_1^+}](v,\tau)+2\mathrm{ch}[{X_1^-}](v,\tau)$.
The intertwiner corresponding to the top gives $o_\cY(v)q^{L_0}$ of the matrix form 
\[\left(\begin{matrix} 0&0&0\\ 0&0&0\\ C&0&0\end{matrix}\right)\left(\begin{matrix} q^n&0&2\pi i q^n\tau\\ 0&q^n&0\\ 0&0&q^n\end{matrix}\right)\]
and so will give a $Z_\cY(v,\tau)$ which is a multiple of $\tau$.

It is highly desirable to understand better the relation with Lyubashenko --- in particular,
why should the SL$(2,\bbZ)$ representations be the same? 

An observation which always holds in the examples is that the graded-dimension $\chi_P(\tau)$ of each projective module
is a  modular form of weight 0 and trivial multiplier for some congruence group. Moreover, it seems in
the examples that the $\bbC$-span of the characters $\chi_P(v,\tau)$ of the projective modules are closed under SL$(2,\bbZ)$.

In the strongly-rational case, it is known that the SL$(2,\bbZ)$-representation contains a principal
congruence group in the kernel. (The principal congruence groups are $\Gamma(N)=\{A\in\mathrm{SL}(2,\bbZ)\,|\,A\equiv I\ (\mathrm{mod}\ N)\}$; modular forms for groups containing
them have by far the richest theory.) In fact, this is exactly what would be expected from the
conjecture of Atkin--Swinnerton-Dyer, which says that if a modular form $f(\tau)$ (with trivial multiplier)  for some subgroup
of SL$(2,\bbZ)$ has a $q$-expansion $f(\tau)=\sum_na_nq^n$ where all coefficients $a_n$ are
algebraic integers, then $f(\tau)$ is a modular form (with trivial multiplier) for a congruence group.  
Surely  the principle congruence groups will continue to play a large role in the strongly-finite theory.

One natural extension of Atkin--Swinnerton-Dyer is as follows. Suppose $\cS$ is a  finite-dimensional (over $\bbC$) space of functions $f(\tau)=q^r\sum_n\sum_if_{n,i}\tau^iq^n$ which is closed under the action of
 SL$(2,\bbZ)$ (at weight-0 say). Suppose in addition that $\cS$ has a basis consisting of 
 functions whose coefficients $f_{n,i}$ are all algebraic integers. Then there are finitely many
 functions $f^{(j)}(\tau)$ in $\cS$, which are modular forms for some $\Gamma(N)$ of some
 nonpositive  weights $k^{(j)}\in\bbZ_{\le 0}$, and the modular closure of these $f^{(j)}(\tau)$ equals $\cS$.
(These $f^{(j)}(\tau)$ would have pure $q$-expansions, but  polynomials in $\tau$ would arise
through the $(c\tau+d)^{-k^{(j)}}$ factors.) This is easily seen to be satisfied by all known examples
of strongly-finite VOAs.

The hope is that modularity arguments can help prove rationality, say, if $\cV$ is already known to
be $C_2$-cofinite.  This would be very useful, because  $C_2$-cofiniteness is easier to show than rationality. It seems reasonable to expect that if the $\bbC$-span of the graded-dimensions $\chi_{M_j}(\tau)$ of
all simple modules of a strongly-finite VOA $\cV$ is closed under SL$(2,\bbZ)$, then $\cV$ is
in fact strongly-rational. In fact, it could be true that if $\cV$ is strongly-finite and its vacuum graded-dimension $\chi_\cV(\tau)$
is fixed by some principal congruence group, then $\cV$ is strongly-rational. The reason this may
hold is because from previous remarks it makes one suspect that $\cV$ is then projective
(as well as being the tensor unit), but as we'll explain below the tensor product of a projective with
any module should be projective, so for any module $M$, $M\cong M\otimes\cV$ will be
projective, and this quickly implies that $\cV$ is strongly-rational. 

Perhaps related is the  interesting conjecture of Gaberdiel--Goddard \cite{GG}, that a strongly-finite VOA is strongly-rational iff 
Zhu's algebra $A_0(\cV)$ is semi-simple. 

\subsubsection{A log-modular tensor category}\label{sec:logmod}

Most fundamental is the question as to what replaces the role of modular tensor category for strongly-finite VOAs.

\begin{defn}\label{def:mod} {A \textbf{log-modular tensor category} $\cC$ is
 a finite tensor category which is ribbon, whose double is isomorphic to the Deligne product $\cC\otimes \cC^{opp}$.}
\end{defn}
 
 These terms are defined in Section \ref{sec:btc}. Modular tensor categories $\cC$ are precisely the ribbon fusion categories
 whose double is isomorphic to $\cC\otimes\cC^{opp}$, and a finite tensor category is the natural nonsemi-simple generalization of
 a fusion category, so this is the natural definition.  A log-modular tensor category
 is meant to be the category of $\cV$-modules for some strongly-finite VOA, so the definition is conjectural in that
 sense and should be tweaked if necessary as we learn more.

There should be several superficially different
but equivalent definitions; this definition should make connection with the deep work of Lyubashenko \cite{L,L2}, and also with the centre construction.

\begin{conjecture}\label{conj:mod}  {The category $\mathrm{Mod}^{g.r}(\cV)$ of $\cV$-modules is a log-modular category.}
\end{conjecture}

The centre (or quantum double) construction applied to a finite tensor category should yield a
    log-modular tensor category (as it does  for modular tensor categories).
Some authors have conjectured that any modular tensor category is realized by a strongly-rational VOA. Analogously,
one should ask if any log-modular tensor category  is realized by a strongly-finite VOA (we don't want to elevate this 
to the status of a conjecture yet).

Conjecture \ref{conj:mod} is not merely abstract nonsense --- it has practical consequences even for the most down-to-earth
researchers. For instance, we now list some of its consequences, which we prove in Section 5 of \cite{CG}. 

\begin{corollary} 
 {If $M$ is simple then $\mathrm{End}(M)=\bbC\mathrm{Id}$. If $M$ is indecomposable then all
elements of $\mathrm{rad(End}(M))$ are nilpotent and $\mathrm{End}(M)=\bbC\mathrm{Id}+\mathrm{rad(End}(M))$.
The projective $\cV$-modules form an ideal in $\mathrm{Mod}^{g.r}(\cV)$, closed under tensor products and taking contragredient.}\end{corollary}

\begin{corollary} {Suppose $\cV$ is strongly-finite but not strongly-rational. Then:}
\begin{itemize}
\item[(a)]   {$\cV$ is not projective as a $\cV$-module.}

\item[(b)]  {$\cV$ has infinitely many families of indecomposable $\cV$-modules, each parametrized
  by complex numbers.} 
  
\item[(c)]   {for any grading-restricted weak $\cV$-module $Y$, $P_i\otimes Y\cong Y\otimes P_i\cong\oplus_{j,k}\cN_{i,j}^k\langle
Y:M_j\rangle P_k$;}

\item[(d)]  {the categorical dimension of any projective $\cV$-module $P$ is 0;}

 \item[(e)] For any projective module $P$ and any simple but non-projective module $A$, the open Hopf link invariant vanishes: $\Phi_{P, A} =0$;
 
\item [(f)] The open Hopf link $\Phi_{P, P'}$ is nilpotent for $P, P'$ projective when $P'$ is indecomposable but not simple. 
\end{itemize}
\end{corollary}

In Section \ref{sec:btc} we define the `semi-simplification' of a category, by quotienting by the negligible morphisms and
restricting to the  full subcategory generated by the simple modules.

\begin{conjecture}  {The semi-simplification of category $\mathrm{Mod}^{g.r}(\cV)$  is a modular tensor category,
whose Grothendieck ring is the quotient of the Grothendieck ring of $\mathrm{Mod}^{g.r}(\cV)$ by the image of the
$\cV$-modules which are both simple and projective.}
\end{conjecture}

This conjecture requires that, when the tensor product of simples in Mod$^{g.r}(\cV)$ is decomposed
into a direct sum of indecomposables, any nonsimple indecomposable arising will have quantum-dimension 0. In the case of the triplet algebras and the even part of the symplectic
fermions, the nonsimple indecomposables are in fact projective. Perhaps that stronger statement continues to hold in general.

To help make this conjecture more concrete,
consider a simple $W$ in any braided finite tensor category (not necessarily satisfying the conjecture). Then  $U\mapsto\Phi_{U,W}$ will be a one-dimensional representation of the tensor ring, 
\begin{equation}\label{fusionrep}
\sum_X {\cN_{U,V}}^X   \Phi_{X,W}=\Phi_{U, W}\circ \Phi_{V, W}\,,
\end{equation}
where the $ {\cN_{U,V}}^X$ are defined as usual by $U\otimes V=\oplus_{X} {\cN_{UV}}^X X$,
and the sum is over the (equivalence classes of) indecomposables $X$ in $\cC$.
 Recall the Hopf link invariant \eqref{hopflink}:
$S^\hopflink_{V,W}=  \mathrm{tr}_W\left( \Phi_{V, W} \right) =S^\hopflink_{W,V}\in \bbC\,.$
Here we take the ordinary (quantum) trace. Since $\Phi_{V,W}$ is a number, $S^\hopflink_{V,W}=
S^\hopflink_{\one,W}\Phi_{V,W}$.
Recall also the  category $\overline{\cC}$ defined at the end of Section 2.6, and
write $\overline{U}$ for the image of $U\in\cC$ in $\overline{\cC}$.
Now, $S^\hopflink_{V,W}=0$ if either the quantum-dimension $S^\hopflink_{\one,V}$ or $S^\hopflink_{\one,W}$ vanishes, in which case $\overline{W}=0$.
 Then for $\overline{W}\in\overline{\cC}$, \eqref{fusionrep} becomes
\[
\sum_X {\cN_{U,V}}^X   \frac{S^\hopflink_{X,W}}{S^\hopflink_{\one, W}}= \frac{S^\hopflink_{U,W}}{S^\hopflink_{\one, W}} \frac{S^\hopflink_{V,W}}{S^\hopflink_{\one, W}}
\]
where the sum is now over all (equivalence classes of) indecomposable objects $\overline{X}$ in $\overline{\cC}$, or equivalently all  equivalence classes of indecomposable objects $X$ in $\cC$
with nonvanishing quantum-dimension. This equation is true for any $\cC$. But if the matrix
 $S^\infty$ is invertible in $\overline{\cC}$ then one obtains the standard Verlinde formula
\begin{equation}\label{eq:standardverlinde}
 {\cN_{U,V}}^X = \sum_{W\in \overline{\cC}} \frac{S^\hopflink_{U,W}S^\hopflink_{V,W}\left(S^{\hopflink-1}\right)_{W,X}}{S^\hopflink_{\one, W}}\,.
\end{equation}
But for this to happen, $S^\hopflink$ would have to be square. So we'd have to restrict $U,V$ in
\eqref{fusionrep}  to simples in $\cC$, in which case any $X$ appearing on the left-side
with nonzero multiplicity would have to be either simple or have quantum-dimension 0.
This is the situation of the previous conjecture, in which case $S^\hopflink$ would coincide
with the modular matrix $S$ of the semi-simplification.

\begin{conjecture}\label{unimodular}  {Any strongly-finite VOA $\cV$ is unimodular, in the sense that the projective cover $P_{{\cV}}$
of the $\cV$-module $\cV$ contains $\cV$ as a submodule (and not merely as a quotient).}
\end{conjecture}
This conjecture is true for log-modular tensor categories according to Proposition 4.5 of \cite{ENO}.

\subsubsection{Open Hopf links and Verlinde}\label{sec:hopfverl}

Most exciting is the question as to what replaces Verlinde's formula. One complication is that there 
are (at least) two natural rings to consider: the Grothendieck ring Fus$^{gr}(\cV)$ and the
simple-projective ring  Fus$^{simp}(\cV)$, discussed in Section \ref{sec:basics}. Verlinde is interested in its regular representation. This can be realised by 
 the fusion matrices $\cN^{gr}_i=(\cN^{gr\,k}_{i,j})$, whose entries are the structure
constants $[[M_i]][[M_j]]=\sum_k\cN^{gr\,k}_{i,j}[[M_k]]$ in Fus$^{gr}(\cV)$. If Fus$^{simp}(\cV)$
is finite-dimensional over $\bbZ$ (and we expect it always is), then the analogous matrices
$\cN^{simp}_M$ can be defined for all indecomposables $M$ in the ring. These fusion matrices
define a representation of Fus$^{gr}(\cV)$ and Fus$^{simp}(\cV)$, namely the regular representation.

As explained in Section 2.4, statement (V1)$'$  is trivial: there is some matrix $S^{gr\,\otimes}$. The first big question is (V2)$'$: interpreting the indecomposable Fus$^{gr}(\cV)$-representations $B^{gr\,h}$.
 In the strongly-rational
case, they are all 1-dimensional, and are given by ratios of Hopf link invariants as explained
in Section 2.3. But in the strongly-finite
case, they are higher dimensional and a new idea is needed. That idea, we propose, is the open Hopf link
invariant.

Let $\mathcal C$ be a 
ribbon category (see Section \ref{sec:btc} for the detailed definitions). 
Of course we are interested in $\cC$ being the category {Mod}$^{g.r}(\cV)$. For each object$=$module $W$, define the map
\[
\Phi_W:\mathcal C \rightarrow \mathrm{End}(W)\,, \qquad U\mapsto \Phi_{U, W}\,,
\]
where $\Phi_{U, W}$ is the \textit{open Hopf link operator} \eqref{openHopf}.
The diagram \eqref{openHopf} translates to the formula 
\begin{equation}\nonumber
\Phi_{V, W} := 
 \mathrm{ptr}_L\left(c_{V, W} \circ  c_{W, V} \right)\in\mathrm{End}(W)\,,
\end{equation}
where the \textit{left partial trace} ptr$_L:\mathrm{End}(V\otimes W)\rightarrow \mathrm{End}(W)$ is
\[
\text{ptr}_L\left(f\right):= 
({d_V}\otimes\Id_W)\circ
(\Id_{V^*}\otimes f)\circ
(b'_{V}\otimes \Id_W)\,,
\]
using maps defined in Section  \ref{sec:btc}. In fact, it isn't hard to show \cite{CG} that each $\Phi_{U,W}$
lies in the centre of End$(W)$.

  A fundamental property of $\Phi_{U,W}$ is that $U\mapsto \Phi_{U,W}$ defines a representation of the tensor ring in the finite-dimensional algebra End$(W)$: $\Phi_{U\otimes V,W}=\Phi_{U,W}\circ\Phi_{V,W}$. The graphical calculus proof of this statement is:
\begin{center}
\begin{tikzpicture}

\node[text width=3cm] at (-2.5, 2.6)    {$\Phi_{U\otimes V, W}$};
\node[text width=3cm] at (-.6, 2.6)    {$=$};

\draw [ar](.38,1) -- (.38, 2.1);
\draw [ar](.3, 3) arc (100:435:.4cm);
\draw (.38,2.3) -- (.38,4);

\node[text width=3cm] at (0, 2.6)    {$U\otimes V$};
\node[text width=3cm] at (1.7, .75)    {$W$};
\node[text width=3cm] at (2.5, 2.6)    {$=$};

\draw [ar] (2.48,1) -- (2.48, 1.55);
\draw (2.48,1.75) -- (2.48, 2.1);
\draw [ar](2.4, 3) arc (100:435:.4cm);
\draw (2.48,2.3) -- (2.48,4);
\draw [ar] (2.38, 3.45) arc (95:440:.9cm);

\node[text width=3cm] at (3.25, 3.6)    {$U$};
\node[text width=3cm] at (3.25, 2.6)    {$V$};
\node[text width=3cm] at (3.8, .75)    {$W$};
\node[text width=3cm] at (5.2, 2.6)    {$=$};

\draw [ar](5.28,1) -- (5.28,1.4);
\draw [ar](5.2, 2.3) arc (100:435:.4cm);
\draw (5.28,1.6) -- (5.28,2.5);

\draw [ar](5.2, 3.4) arc (100:435:.4cm);
\draw (5.28,2.7) -- (5.28,4);

\node[text width=3cm] at (5.8, 3)    {$U$};
\node[text width=3cm] at (5.8, 1.9)    {$V$};
\node[text width=3cm] at (6.6, .75)    {$W$};

\node[text width=3cm] at (7.5, 2.6)    {$=$};
\node[text width=3cm] at (8.2, 2.6)    {$\Phi_{U, W}\circ\Phi_{V, W}\,$.};

\end{tikzpicture}\end{center}
Here, the third equality holds because of naturality of braiding, and the second one is \eqref{braiding} together with $\mathrm{Id}_{U\otimes W}=\mathrm{Id}_U\otimes \mathrm{Id}_W$.
These $\Phi_{U,W}$ obey several other properties. For example, if $\cC$ has direct sums
(and ours always will), then 
$\Phi_{V,W_1\oplus W_2}=\Phi_{U,W_1}\oplus\Phi_{V,W_2}$. In other words,
 without loss of generality we can restrict to indecomposable $W$.

In some cases at least, this representation $U\mapsto\Phi_{U,W}$ of Fus$^{full}(\cV)$
descends to a representation of Fus$^{gr}(\cV)$ (at present we don't understand when this happens).
We propose that (V2)$'$ is:

 \begin{conjecture}  
  Each indecomposable  subrepresentation $M\mapsto B^{h}_M$ of the
 Grothendieck
 ring representation $\cN^{gr}_M$ is equivalent to a subrepresentation of the open Hopf link representation $M\mapsto \Phi_{M,P_j}$ for some  projective cover $P_j$ of a simple module $M_j$.  Each indecomposable subrepresentation of the
 regular representation of the Grothendieck ring Fus$^{gr}(\cV)$ or the full tensor ring Fus$^{full}(\cV)$ occurs  with multiplicity 1. \end{conjecture}

This is the complete story for (V2)$'$ for the Grothendieck ring, except we don't know 
which Fus$^{full}(\cV)$-representations $\Phi_{\star,W}$ descend to Fus$^{gr}(\cV)$
(perhaps all do).
 We do not know which indecomposable subrepresentations should appear in the regular representation of the full tensor ring Fus$^{full}(\cV)$ or even the 
 subring Fus$^{simp}(\cV)$ generated by the simple $\cV$-modules $M_j$ and their projective
 covers $P_j$ --- even for
 $\cV=\cW_p$ this question is  mysterious. It is a theorem for $\cW_p$ that open Hopf link operators do \textit{not} suffice
 for  Fus$^{full}(\cV)$ or  Fus$^{simp}(\cV)$. This interesting question should be explored. 
The open Hopf link is a vertical arrow through the trivial knot; could the operator-valued invariants associated to arrows through certain nontrivial knots be responsible for the missing indecomposables?

Now turn to (V3)$'$, which is the heart of Verlinde. For the $\cW_p$ models, \cite{PRR} relate 
$S^{gr\,\otimes}$ to 
the $-1/\tau$ transformation of (pseudo-)characters --- we describe this interesting observation later this section,
but it isn't clear at present how it generalizes to other strongly-finite $\cV$. In any case, at present
we have nothing to add to the relation of  $S^{gr\,\otimes}$ to the modular $S$-matrix, although those $\cW_p$
calculations suggest something general can be said.

At least as  subtle  is to relate the open Hopf link operators to the modular $S$-matrix. Because
the latter is numerical (at least once a basis of 1-point functions is chosen), it is probably wise to
see how to extract numerical invariants from the endomorphisms $\Phi_{U,W}$. 
 
Let us sketch one way which underlies one of our observations. 
Let  $\mathcal P$ be the (full) subcategory of $\cC$ consisting of the projective objects in $\cC$. Then among other things, $\cP$ is semi-simple (with irreducibles being the projective covers of
the simple objects of $\cC$), and $\cP$ is a tensor ideal of $\cC$ (since the tensor product of a
projective with any object of $\cC$ will be projective). Suppose  $\cC$ contains a projective simple
self-dual object $\tilde{P}$. Then, provided $\cC$ is unimodular in the sense of
Conjecture \ref{unimodular}, there will be a \textit{modified trace}  \cite{CGP,GKP1} on $\cP$: i.e.  
 a family of linear functions 
\[
\left\{\ t_V:\mathrm{End}(V) \rightarrow \mathbb C\ |\ V\in\mathcal P\ \right\}
\]
satisfying
\begin{enumerate}
\item For $U$ in $\mathcal P$ and $W$ in $\mathcal C$ and any $f$ in $\mathrm{End}\left(U\otimes W\right)$ 
\[
t_{U\otimes W}\left(f\right) =t_W\left( \mathrm{ptr}_L\left(f\right)\right)\,.
\]
\item For $U, V$ in $\mathcal P$ and morphisms $f:V\rightarrow U$ and $g:U\rightarrow V$
\[
t_V\left(g\circ f\right)=t_U\left(f\circ g\right)\,. 
\]
\end{enumerate}
This idea of modified trace has been developed by Geer, Patureau-Mirand and collaborators \cite{GKP1, CGP, GKP2, GPT, GPV}.
 
Given a modified trace on $\cP$, we can define the {\it logarithmic Hopf link
invariant}
\[
S^{\hopflink\,;\cP,x}_{V, W}:= t_W\left(\Phi_{V, W} \circ x \right)=t_W\left(x\circ \Phi_{V,W}\right)\in
\bbC\,,
\]
for any $V$ in $\mathcal C$, $W$ in $\mathcal P$, and $x$ in $\mathrm{End}\left(W\right)$.
For readability, we will drop the endomorphism $x$ from the superscript if $x$ is the identity, and we will
drop the $\cP$ if $x$ is present.
The logarithmic Hopf link invariant is symmetric just like the non-logarithmic one: 
$S^{\hopflink\,;\cP}_{V, W} =S^{\hopflink\,;\cP}_{W, V}$. Moreover, for any $x\in \textrm{End}(W)$,
\begin{equation}\label{Philog}
S_{U, W}^{\hopflink\,; x}=t_W \left(\Phi_{U, W}\circ x \right) =  t_W \left( \Phi_{\one, W}\circ(\Phi_{U, W}\circ x)\right)
= S_{\one, W}^{\hopflink\,;\Phi_{U, W}\circ x}\,. 
\end{equation}

Again, logarithmic Hopf links are a promising way to extract numbers from  the endomorphisms
$\Phi_{U,W}$.
More generally, if $\cP$ is any tensor ideal of $\cC$ and $\tilde{P}$ is a right ambidextrous object,
then we will get a modified trace on $\cP$ and hence can define logarithmic Hopf links.

Incidentally, there is a relation between the pseudo-traces used by Miyamoto, and these
modified traces, and we expect this to play a role in the general story. Let $\cA_\cV$ be the
finite-dimensional algebra of section \ref{sec:reptype} (see also Theorem 2.3 of \cite{CG}). Recall that given any symmetric linear functional $\phi$
on $\cA_\cV$, we get a pseudo-trace $\phi_P :\mathrm{End}_{\cA_\cV}(P)\rightarrow\bbC$
 in the usual way
for each projective finite-dimensional $\cA_\cV$-module
$P$. 
{Given any modified trace $\{t_P\}$, there is one and only one symmetric linear functional $\phi$ 
on $\cA_{\cV}$,
 such that $t_P=\phi_P$ for all projective $\cV$-modules $P$. Conversely, the pseudo-traces $\phi_P$ of any
symmetric linear functional $\phi$ on $\cA_\cV$ necessarily satisfy property 2 of the definition of modified trace (but not in general property 1).} See Section 6.1 of \cite{CG} for more details.

There is a striking relation between the modular matrix $S^\chi$ and the logarithmic Hopf
link invariants, for the $\cW_p$ models, and we expect this to generalize to any strongly-finite
$\cV$. A basis for the $3p-1$-dimensional space of 1-point functions at weight 0 is $\mathrm{ch}[P^+_\ell],\mathrm{ch}[X^+_\ell], \mathrm{pch}[X^+_\ell]$ for all $\ell=1,2,\ldots,p-1$, together with
$\mathrm{ch}[X^\pm_p]$ for both signs. By $S^\chi$ we mean the matrix corresponding
to the modular transformation $\tau\mapsto -1/\tau$ in this basis. Then the entries of $S^\chi$
appropriately normalized agree with logarithmic Hopf link invariants, appropriately normalized.
We make this explicit shortly (see equations (\ref{S/S=S/S}, \ref{cS=cS})), but two typical examples are
\begin{equation}\label{ex:S}
\frac{S^\chi_{P^+_s, P^+_\ell}}{S^\chi_{X^+_1, P^+_\ell}} = \frac{S^{\hopflink\,;\cP}_{P^+_s, P_\ell^+}}{S^{\hopflink\,;\cP}_{X^+_1, P^+_\ell}}\,,\qquad \frac{S^\chi_{X^0_s, P^+_\ell}}{S^\chi_{X^0_1, P^+_\ell}} = \frac{S^{\hopflink\,;x}_{X^+_s, P_\ell^+}}{S^{\hopflink\,;x}_{X^+_1, P^+_\ell}}
\end{equation}
where $X^0_s$ refers to pch$[X^+_\ell]$.

This interpretation of  $S^\chi$ in terms of $S^\hopflink$, combined with the expression due to \cite{PRR} of $S^\otimes$ in terms of $S^\chi$, allows us to write $S^\otimes$ in terms of
the (logarithmic) Hopf link invariants. 
And our  (V2)$'$ conjecture tells us we can interpret the indecomposable subrepresentations
$B^h$ of the regular representation of Fus$^{gr}(\cV)$ in terms of $\Phi_{U,W}$, i.e. in terms of
the (logarithmic) Hopf link invariants. Putting these observations
together, we get a Verlinde formula akin to \eqref{VerlForm}, for the Grothendieck ring.

There are other  ways to arrive at Verlinde-like formulas. One was discussed in Section 3.1.2.
For another example, assume that a projective module $W$ has 2-dimensional endomorphism ring 
(examples are $P^\pm_s$ of $\cW_p$). Then it almost has a preferred basis:  $\text{Id}_W$ and 
some $x=x_W$ with $x^2=0$. Then every open Hopf link operator is of the form
\[
\Phi_{U, W}=a_{U, W}\mathrm{Id}_W  +b_{U, W} x_W
\]
for complex numbers $a_{U, W}$ and $b_{U, W}$. These numbers relate to the tensor product structure constants via
\begin{equation}
\begin{split}\label{eq:verlinde}
\sum_X {\cN_{U,V}}^X a_{X, W} = a_{U, W}a_{V, W}, \qquad
\sum_X {\cN_{U,V}}^X b_{X, W} = a_{U, W}b_{V, W}+ b_{U, W}a_{V, W}
\end{split}
\end{equation}
where the sums are over all indecomposable objects. 
On the other hand these numbers also relate to the logarithmic Hopf link invariants, since
\[
S^{\hopflink\,;\cP}_{U, W}=a_{U,  W}S^{\hopflink\,;\cP}_{\one, W}+b_{U, W} S^{\hopflink\,; x}_{\one, W}\,, \qquad S^{\hopflink\,; x}_{U, W}=a_{U,  W}S^{\hopflink\,;x}_{\one, W}\,,
\]
so that they can be expressed in terms of normalized logarithmic Hopf link invariants:
\begin{equation}\nonumber
\begin{split}\label{eq:qdims}
a_{U, W}=\frac{S^{\hopflink\,; x}_{U, W}}{S^{\hopflink\,; x}_{\one, W}}\,, \qquad
b_{U, W}=  \frac{S^{\hopflink\,;\cP}_{\one, W}}{S^{\hopflink\,; x}_{\one, W}} \left(\frac{S^\hopflink_{U, W}}{S^{\hopflink\,;\cP}_{\one, W}}- \frac{S^{{\hopflink\,; x}}_{U, W}}{S^{\hopflink\,;x}_{\one, W}}\right),
\end{split}
\end{equation} 
assuming the numbers $S^{\hopflink\,; x}_{\one, W}$ don't vanish. One could then hope to
invert \eqref{eq:verlinde} for some tensor subcategory of $\cC$, expressing ${\cN_{U,V}}^X$
in terms of the logarithmic Hopf link invariants. For instance, in the triplet example, one recovers
the standard Verlinde formula \eqref{eq:standardverlinde}, for the subcategory generated
 by the simples.   In that example we also get a Verlinde-type formula for the fusion coefficients of type simple with projective, see \eqref{eq:vertrip} below.

\vspace{5mm}

\subsection{The picture in the locally finite setting}

An abelian category is called \textit{locally finite} if Hom$(X, Y)$ is finite-dimensional for any two objects in the category and if all objects have finite Jordan-H\"older length. 
Logarithmic CFTs with locally finite module categories are quite common and the recent development here has helped lead to our perspective. Older developments are already reviewed in the introduction to logarithmic CFT \cite{CR4} and the reviews \cite{RW2, C}.

 The starting point had been logarithmic conformal field theories with a continuum of simple objects. Then characters are not modular in the sense that they are not components of a finite-dimensional vector-valued modular form. Nonetheless it is still possible to express the $S$-transformation $\tau \mapsto -1/\tau$ as an integral of characters:
 \[
 \text{ch}[X](-1/\tau) = \int_I S_{X, Y} \text{ch}[Y](\tau) dY. 
 \]
Here $I$ is the index set parameterizing simple objects. The $S$-matrix is thus replaced by an $S$-kernel $S_{X, Y}$.
What should then replace the Verlinde formula? It is a natural guess that this would be
\[
N_{X Y }^{\ \ \ Z} = \int_I \frac{S_{X, W}S_{Y, W}S^{-1}_{W, Z}}{S_{\cV, W}} dW.
\]
$\cV$ is the VOA itself. Here of course one needs to make sense of the inverse of the $S$-kernel. This is to be interpreted in a distributional sense and so are the fusion coefficients $N_{X Y }^{\ \ \ Z} $. Further these formulas are subtle in the sense that one has to choose the integration measure wisely. However it turns out that in many cases, that is the free boson \cite{CR4}, the affine Lie superalgebra of $\mathfrak{gl}(1|1)$ \cite{CQS1, CR3}, its extensions \cite{AC} and especially in the fractional level WZW theories of $\mathfrak{sl}_2$ \cite{CR1, CR2}, these ideas work extremely nice. Modular objects involved here include meromorphic Jacobi forms, distributions and mock modular forms. 

The next natural question is now whether there is a good explanation for these experimental observations. In \cite{CM1} the singlet algebra was studied, see also \cite{RW1}. This is a $U(1)$-orbifold of the triplet and conversely the triplet is an infinite order simple current extension of the singlet. 
Here characters are false theta functions and there modularity had not even been known. A regularization, called $\epsilon$, was needed to obtain such modular properties. It then turned out that regularized asymptotic dimensions of characters
\[
\text{qdim}[X, \epsilon] = \lim_{\tau\rightarrow 0} \frac{\text{ch}[X, \epsilon](\tau)}{\text{ch}[\cV, \epsilon](\tau)}
\]
capture the tensor ring. Moreover, depending on the regime of the regularization parameter it actually either captures the semi-simplification or the full Grothendieck ring. This has then been further and deeper investigated in \cite{CMW} and generalized to higher rank analogues \cite{CM2}. 

The for the moment final natural question is then to give a categorical meaning to these findings. The answer is: Asymptotic dimensions in the continuous regularization parameter regime correspond to logarithmic Hopf links while the regime that captures the semi-simplification corresponds to logarithmic Hopf links with nilpotent endomorphism insertion \cite{CMR}. 
This is to be taken as a sign that the interplay between modular-like tensor categories, modular forms and logarithmic CFT goes even beyond the strongly-finite setting and some kind of generalization to the locally-finite setting will be possible.

\subsection{Example: the triplet algebra $\cW_p$}\label{sec:V}

In this subsection we verify explicitly many of our conjectures  for the $\cW_p$ algebras.
The other class of examples which can be studied this explicitly are the even part $SF^+_d$
of the symplectic fermions. These two families nicely complement each other: $\cW_p$ has
an unbounded number of simples but $L_0$ has nilpotent order 2, while $SF^+_d$ always  has
4 simples but its $L_0$ has arbitrarily high nilpotent order.

In most ways the symplectic fermions are simpler. But
one place they differ is that the modular closure of the $2p$-dimensional space spanned
by the $\cW_p$-characters recovers the full $3p-1$-dimensional SL$(2,\bbZ)$-module $\cF(\cW_p)$
expected from both the Miyamoto and Lyubashenko pictures. By contrast, the modular closure of
the 4-dimensional space spanned by the $Sp^+_d$-characters is  $d+4$-dimensional,
while the dimension of $\cF(SF^+_d)$ is expected to be $2^{2d-1}+3$.
Another difference, as we'll see in Section 3.6, is that $\cW_p$ is tame while $SF^+_d$, like
most strongly-finite VOAs,  has wild-type.

\subsubsection{Calculating the open Hopf link invariants}

Crucial to our story are the open Hopf link invariants. These are difficult to compute in Mod$^{g.r}(\cV)$
because we don't yet have a good grasp on the braiding. However,
\cite{FGST2} originally conjectured that the representation categories of $\cW_p$ and  the restricted quantum group $\Uq$ at a $2p$-th root of unity are equivalent as braided tensor categories.
\cite{NT} proved equivalence as abelian categories (the case $p=2$ was also done in \cite{FGST2}), and it is a direct comparison of  \cite{TW} with  \cite{KS}  that the corresponding tensor product rings 
Fus$^{simpl}$ of simples and projectives are isomorphic. Braided equivalence is however false \cite{KS}, but it is believed that the difference is minor, only lieing in different associativity isomorphisms. This is 
proven in the case $p=2$ by Gainutdinov and Runkel \cite{GaR} using \cite{DR}. Although these 
braided tensor categories are not equivalent, it is expected that the
open Hopf link invariants of greatest concern to us will be the same.

The quantum group $\Uq$ is described for example in  Murakami \cite{Mu}. Original references are \cite{FGST1,FGST2, KS, S}.
Let $q=e^{\pi \Ii/p}$ be a primitive $2p$-th root of unity, and write
$\{n\}_q:= q^n-q^{-n}=2i\sin(\pi/p)$.
The restricted quantum group $\Uq$ doesn't have a universal $R$-matrix, i.e. it is not a quasi-triangular Hopf algebra. However it can be embedded in a quasi-triangular Hopf algebra $\overline D$  with universal $R$-matrix \cite{KS,FGST1}. All 
modules of $\Uq$ that we are interested in, lift to $\overline D$, and
the universal $R$-matrix action of $\overline D$ on them coincides with that of $\Uq$.

The simple $\Uq$ modules are $U^\pm_s$ for either sign and $1\le s\le p$. Their projective
covers are $R^\pm_s$, where $R^\pm_p=U_p^\pm$ (in the literature, the $R^\pm_s$ are often denoted by $P^\pm_s$, but we chose  `$R$'  to avoid confusion with $\cW_p$). The correspondence
between the $\Uq$- and $\cW_p$-modules of interest to us are: $U^\pm_s\leftrightarrow X^\pm_s$
and $R^\pm_s\leftrightarrow P^\pm_s$. For example, the Loewy diagrams of $R^\pm_s$ and $P^\pm_s$ are
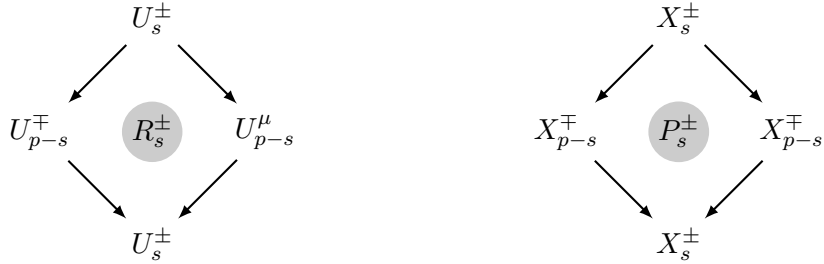
\begin{center}
\begin{tikzpicture}[thick,>=latex,nom/.style={circle,draw=black!20,fill=black!20,inner sep=1pt}]
\node (top1) at (5,1.5) [] {$U^\pm_s$};
\node (left1) at (3.5,0) [] {$U^\mp_{p-s}$};
\node (right1) at (6.5,0) [] {$U^\mu_{p-s}$};
\node (bot1) at (5,-1.5) [] {$U^\pm_s$};
\node at (5,0) [nom] {$R^\pm_s$};
\draw [->] (top1) -- (left1);
\draw [->] (top1) -- (right1);
\draw [->] (left1) -- (bot1);
\draw [->] (right1) -- (bot1);

\node (top2) at (12,1.5) [] {$X^\pm_s$};
\node (left2) at (10.5,0) [] {$X^\mp_{p-s}$};
\node (right2) at (13.5,0) [] {$X^\mp_{p-s}$};
\node (bot2) at (12,-1.5) [] {$X^\pm_s$};
\node at (12,0) [nom] {$P^\pm_s$};
\draw [->] (top2) -- (left2);
\draw [->] (top2) -- (right2);
\draw [->] (left2) -- (bot2);
\draw [->] (right2) -- (bot2);
\end{tikzpicture}
\end{center}
The endomorphism ring of each $U^\pm_s$, $1\le s\le p$, can be canonically identified with $\bbC$.
The endomorphism ring of each projective module $R^\pm_s$, $1\le s<p$, is spanned by the identity 
(which we'll denote by $e$) as well as a nilpotent one we'll call $x$,  whose image is the simple socle module of its Loewy diagram. 

The open Hopf link operators $\Phi_{\star,U^\pm_s}$ for  $1\le s<p$ are easy to compute,
and we find
\[\Phi_{U^{\epsilon'}_{s'}, U^\epsilon_s}=(-1)^{s'+1}\epsilon''\frac{\sin(\pi ss'/p)}{\sin(\pi s/p)}\,,\qquad 
\Phi_{R^{\epsilon'}_i, U^\epsilon_s}=0\,,\]
for all  $1\le s<p$,  $1\le s'\le p$, $1\le i<p$, and $\epsilon,\epsilon'\in\{\pm 1\}$, where $\epsilon''=+1$ unless
$\epsilon=\epsilon'=-$ in which case $\epsilon''=-1$. From this we obtain
\[S^\hopflink_{U^{\epsilon'}_{s'}, U^\epsilon_s}=S^\hopflink_{U^{\epsilon}_{s}, U^{\epsilon'}_{s'}}=(-1)^{s'+1}\epsilon\epsilon'\frac{\sin(\pi ss'/p)}{\sin(\pi /p)}\,,\qquad
S^\hopflink_{R^{\epsilon'}_{i}, U^\epsilon_s}=S^\hopflink_{ U^\epsilon_s,R^{\epsilon'}_{i}}=0\,.\]
In particular, $S^{ss}_{U^{\epsilon'}_{s'}, U^\epsilon_s}=S^{aff}_{U^{\epsilon'}_{p-s'}, U^{\epsilon}_{p-s}}$.

Incidentally, the semi-simplification yields the fusion ring of the modular tensor category corresponding to $L_{p-2}(\mathfrak{sl}_2)\otimes L_1(\mathfrak{sl}_2)$, but the semi-simplification Mod$^{g.r}(\cW_p)^{ss}$ itself
cannot be identified with that of the affine algebra, but rather with its \textit{conformal flow}  in the sense of \cite{LM}, and in particular is a nonunitary modular tensor category, as expected.

The open Hopf links for the projectives is much more subtle. In \cite{CG} we compute them
using two rather different methods. One method followed  \cite{CGP} and uses the known tensor ring, while the other one is via a deformation argument \`a la Murakami--Nagatomo \cite{MN, Mu} and it does not assume knowledge of the tensor ring. This is important if one takes our point of view that logarithmic Hopf link invariants essentially determine the tensor ring. The idea of deformation seems to be quite fruitful: the logarithmic Hopf links of the unrolled restricted quantum group $\overline{\text{U}}_q^H(\mathfrak{sl}_2)$ can be computed using the idea of deformable families of quantum group modules \cite{CMR}, while much of the structure of coset vertex algebras can be deduced by treating them as deformable families of VOAs \cite{CL, CL2}. In the case of $L_k\left(\mathfrak{gl}\left(1|1\right)\right)$ it is also possible to construct deformable families of modules, and the usefulness of this observation should be further explored.

We find that the remaining open Hopf link invariants are \cite{CG}:
\begin{theorem}\label{openhopfWp}
The open Hopf link operators of  $R_j=R_j^+$ and $R^-_{p-j}$  are the same.
For any $1\le s\le p$, $1\le i,j<p$, and any signs $\epsilon,\epsilon'$,
\begin{align}\nonumber
\Phi_{U_s^{\epsilon'}, U_p^\epsilon}&= \epsilon^{j+1}\epsilon''j\,, \qquad\qquad\qquad \Phi_{R_i^{\epsilon'}, U_p^\epsilon}= \epsilon^{i+1}\epsilon''2p\,,\\ \nonumber
\Phi_{U^\epsilon_i, R_j}&= \frac{(-1)^{i+1}\epsilon^{p+j}}{\{j\}_q}\left( \{ij\}_qe+\left((i-1)\{(i+1)j\}_q-(i+1)\{(i-1)j\}_q\right)x\right)\,,  \\
\nonumber \Phi_{R^\epsilon_i, R_j}&=(-1)^{i+1}\epsilon^{p+j}4p \cos(\pi ij/p)\,x\,,
\end{align}
where $\epsilon''=1$ if $\epsilon=+$ otherwise $\epsilon''=\epsilon^{\prime\,p}$.
\end{theorem}

Likewise, the logarithmic Hopf link invariants of $R_i=R^+_i$ and $R^-_{p-i}$ coincide. We obtain: 
\begin{theorem}\label{thm: logHopf} For any $1\le s\le p$, $1\le i,j<p$, and sign $\epsilon$,
\begin{equation}\nonumber
\begin{split}
\qquad\qquad  
S^{\hopflink\,;\cP}_{U_i^\epsilon, R_j}&= (-1)^{i+j+1}\epsilon^{p+j}2i\cos(\pi ij/p)\,,\qquad
S^{\hopflink\,;x}_{U_i^\epsilon, R_j}= (-1)^{i+j+1}\epsilon^{p+j}\frac{\sin(\pi ij/p)}{\sin(\pi j/p)},
\\
S^{\hopflink\,;\cP}_{R^\epsilon_i, R_j} &= (-1)^{i+j+1}4p\cos(\pi ij/p)\,,\qquad S^{\hopflink\,;x}_{R_i^\epsilon, R_j} =0\,,\qquad
S^{\hopflink\,;\cP}_{U^\epsilon_i, U^+_p}= (-1)^{p-1}i\,,\\  S^{\hopflink\,;\cP}_{R_i^\epsilon, U^+_p}&=(-1)^{p-1}2p\,, \qquad
S^{\hopflink\,;\cP}_{U^\epsilon_i, U^-_p}= (-1)^{i-1}\epsilon^pi\,,\qquad S^{\hopflink\,;\cP}_{R_i^\epsilon, U^-_p}=(-1)^{i-1}\epsilon^p2p\,.\\
\end{split}
\end{equation}
\end{theorem}

\subsubsection{Block-diagonalizations for $\cW_p$}\label{sec:V1}

In this subsection we review the block-diagonalization of the (regular representation of the) Grothendieck ring of $\cW_p$ performed by \cite{PRR},
and  the block-diagonalization of the tensor ring of $\cW_p$ performed by \cite{R}.

Consider first the Grothendieck ring of $\cW_p$. It has basis $X^\pm_s$, which we put in the order 
$X^+_1,X^-_1,X^+_2,X^-_2,\ldots,X^+_p,X^-_p$. As the tensor ring is generated by the simple-current $X^-_1$ and $X^+_2$, it suffices to
block-diagonalize $J^{gr}:=\cN^{gr}_{X^-_1}$ and $Y^{gr}:=\cN^{gr}_{X^+_2}$.

Consider first $p=2$. Then 
\[Q^{gr-1}JQ^{gr}=\mathrm{diag}(1,-1,-1,1)\,,\ Q^{gr-1}YQ^{gr}=\mathrm{diag}(2,\left(\begin{matrix}0&1\\ 0&0\end{matrix}\right),-2)\,,\]
for
\[Q^{gr}=\left(\begin{matrix}1&4&0&-1\\ 1&-4&0&-1\\ 2&0&4&2\\ 2&0&-1&2\end{matrix}\right)\,.\]

More generally \cite{PRR}, there is an invertible matrix $Q^{gr}$ (discussed below) which simultaneously diagonalizes $J^{gr}$ and 
puts $Y^{gr}$ into Jordan canonical form. Those matrices $Q^{gr-1}J^{gr}Q^{gr}$ and $Q^{gr-1}Y^{gr}Q^{gr}$ fall naturally into two $1\times 1$ 
blocks and $p-1$ $2\times 2$ blocks:
\[Q^{gr-1}J^{gr}Q^{gr}=\mathrm{diag}(1;-I_2;\ldots;(-1)^jI_2;\ldots;(-1)^{p-1}I_2;(-1)^p)\,,\]
\[ Q^{gr-1}Y^{gr}Q^{gr}=\mathrm{diag}(\lambda_0;B_{\lambda_1,2};\ldots;
B_{\lambda_j,2};\ldots;B_{\lambda_{p-1},2};\lambda_p)\,,\]
where $1\le j<p$, 
 the eigenvalues are $\lambda_j=2\cos(\pi j/p)$, and we write $I_k$ for the $k\times k$ identity matrix and
 $B_{\lambda,k}$ for the canonical $k\times k$ Jordan block with eigenvalue $\lambda$.

Because any other matrix $\cN^{gr}_{X^\pm_s}$ will be a polynomial in $J^{gr}$ and $Y^{gr}$, $Q^{gr}$ also block-diagonalizes
all $\cN^{gr}_{X^\pm_s}$. The $0$-th and $p$-th blocks  of $Q^{gr-1}\cN^{gr}_{X^\pm_s}Q^{gr}$ will be the numbers $s$ and $(\pm 1)^p(-1)^{s-1}s$ respectively; its $j$-th block, for $1\le j< p$, will be $2\times 2$ upper-triangular,
with diagonal entries (eigenvalues) $\pm \sin(\pi js/p)/\sin(\pi s/p)$. Those $2\times 2$ blocks will not in general
be canonical Jordan blocks, however.

Now turn to the tensor ring Fus$^{simp}(\cW_p)$ of $\cW_p$. 
We are interested in the full subcategory spanned by the irreducible $\cW_p$-modules $X_s^\pm$ and their projective covers
$P_s^\pm$. The corresponding $4p-2$ tensor matrices $\cN_M$, each of size $(4p-2)\times(4p-2)$, were put into simultaneous
block form in \cite{R}. We will put these $4p-2$  $\cW_p$-modules in order $X^+_1,X^-_1,\ldots,X^+_p,
X^-_p,P^+_1,P^-_1,\ldots,P^+_{p-1},$ $P^-_{p-1}$. Again, it suffices to
block-diagonalize $J:=\cN_{X^-_1}$ and $Y:=\cN_{X^+_2}$. 

Let us consider $\cW_2$ first. A matrix block-diagonalizing all 6 tensor matrices is
\[Q=\left(\begin{matrix}1&1&1&0&0&1\cr 1&1&-1&0&0&1\cr 2&0&0&1&0&-2\cr 2&0&0&-1&0&-2\cr 4&0&0&0&1&4\cr 4&0&0&0&-1&4\end{matrix}\right)\,.\]
Then \cite{R}
\[Q^{-1}JQ=\mathrm{diag}(1,1,-1,-1,-1,1)\,,\qquad Q^{-1}YQ=\mathrm{diag}(2,0,B_{0,3},-2)\,.\]

More generally \cite{R}, 
there is an invertible matrix $Q$ (discussed below) which simultaneously diagonalizes $J$ and puts $Y$ into Jordan 
canonical form. Those matrices $Q^{-1}JQ$ and $Q^{-1}YQ$ fall naturally into two $1\times 1$ blocks and $p-1$ pairs of $1\times 1$ and $3\times 3$ blocks:
\[Q^{-1}JQ=\mathrm{diag}(1;1,-I_3;\ldots;(-1)^{j-1},(-1)^jI_3;\ldots;(-1)^{p-2},(-1)^{p-1}I_3;(-1)^p)\,,\]
\[Q^{-1}YQ=\mathrm{diag}(\lambda_0;\lambda_1,B_{\lambda_1,3};\ldots;\lambda_j,B_{\lambda_j,3}; \lambda_{p-1},B_{
\lambda_{p-1},3};\lambda_p)\,.\]
Because  any other matrix $\cN_{X^\pm_s}$ or $\cN_{P^\pm_s}$ is a polynomial in $J$ and $Y$, $Q$ also block-diagonalizes
them. The 0th and $p$th blocks  of $Q^{-1}\cN_{X^\pm_s}Q$ will be the numbers $s$ and $(\pm 1)^p(-1)^{s-1}s$ 
respectively; the $1\times 1$ part of its $j$th block, for $1\le j< p$, will be the number $(\pm 1)^{j-1}\sin(\pi js/p)/\sin(\pi s/p)$; 
the $3\times 3$ part of its $j$th block will be  upper-triangular,
with diagonal entries (eigenvalues) $(\pm 1)^j\sin(\pi js/p)/\sin(\pi s/p)$. 
Again, those $3\times 3$ blocks will not in general be canonical Jordan blocks, however. 
 The block structure for $Q^{-1}\cN_{P^\pm_s}Q$ is similar.

\cite{PRR} found an interesting interpretation of $Q^{gr}$ in terms of the modular $S$-matrix.
Choose the basis ch$[{X^+_s}]$, pch$[X^+_\ell]$, ch$[{X^-_s}]$, for $1\le s\le p$ and $1\le\ell<p$ and
either sign. This is slightly different than the choice we make. Let $\tilde{S}^\chi$ denote the corresponding matrix realizing the modular transformation $\tau\mapsto-1/\tau$.
Then
\begin{equation}
Q^{gr}=\left(\tilde{S}^\chi_{\updownarrow,X^+_s},\,\tilde{S}^\chi_{\updownarrow,X^0_\ell}\,\tilde{S}^\chi_{\updownarrow,\,X^-_\ell},\tilde{S}^\chi_{\updownarrow,X^-_\ell}\right)\end{equation}
where we use the subscript 0 for the pseudo-character.

\subsubsection{Verlinde (V2)$'$ for $\cW_p$}\label{sec:V2}

In this subsection we identify the indecomposable representations contained in the regular representation of
the  Grothendieck and tensor rings Fus$^{gr}(\cW_p)$, Fus$^{simp}(\cW_p)$ of $\cW_p$, with open Hopf link operators. We gave their decomposition last subsection, expressing them as direct sums of 1-, 2- and 3-dimensional
subrepresentations.

The Grothendieck ring and tensor ring have the same two 1-dimensional subrepresentations, corresponding to $j=0$ and 
$j=p$.  The $j=0$ one is uniquely determined by it sending $X^-_1$ to 1 and $X^+_2$ to 2; while the $j=p$ one
sends $X^-_1$ to $(-1)^p$ and $X^+_2$ to $-2$. Using Theorem \ref{openhopfWp}, we identify these with
$\Phi_{\star,U_p^+}$ and $\Phi_{\star,U_p^-}$ respectively.

 The Grothendieck ring also has a 2-dimensional subrepresentation for each $1\le j<p$, which sends $X^-_1$ to 
 $(-1)^jI_2$ and $X^+_2$ to $B_{\lambda_j,2}$. By Theorem \ref{openhopfWp}, 
  this representation is isomorphic to $\Phi_{\star,R^+_{p-j}}$ or equivalently $\Phi_{\star,R^-_j}$.

The tensor ring Fus$^{simp}(\cW_p)$ has a  1-dimensional representation for each $1\le j<p$, sending $X^-_1$ to $(-1)^{j-1}$ and $X^+_2$ to $\lambda_j$. This is  $\Phi_{\star,U^{\epsilon}_{p-j}}$,
where $\epsilon$ is the sign of $(-1)^{j-1}$.

What remains is to understand the 3-dimensional representations in the tensor ring. However:

\begin{theorem} \cite{CG} For any $\cW_p$-modules $U,M$,  the map $\Phi_{U,M}\in\mathrm{End}(M)$ decomposes only into
$1\times 1$ and $2\times 2$ Jordan blocks. In particular, no $\Phi_{\star,M}$ can contain a subrepresentation equivalent to any of the 3-dimensional
indecomposable tensor ring representations given in the previous subsection.\end{theorem}

In Section 3.1.3 we proposed some generalizations of the open Hopf link which may work,
but we have not investigated this further at this point.

Incidentally, the 4-dimensional Grothendieck ring of $SF^+_d$ likewise decomposes into 1- and 2-dimensional indecomposables, which are in turn indecomposable summands in open Hopf link
invariants of simples and projectives. Again, the indecomposables of the regular representation
of $SF^+_d$ are not all subrepresentations of any open Hopf link invariant $\Phi_{\star,W}$.

\subsubsection{Verlinde (V3)$'$ for $\cW_p$}\label{sec:V3}

Comparing the open Hopf link calculations in Section 3.3.1 with the modularity calculations
in ection 2.8, we find that the (logarithmic) Hopf link invariants of the restricted quantum group and 
modular $S$-matrix coefficients for the (pseudo-)characters satisfy the remarkable formula
\begin{equation}\label{S/S=S/S}\frac{S^\chi_{Y,Y'}}{S^\chi_{Y'',Y'}}=\frac{S^{\hopflink\,;\gamma}_{Y_U,Y'_U}}{S^{\hopflink\,;\gamma}_{{Y_U''},{Y_U'}}}\end{equation}
for any $Y,Y'\in\{P^+_s,X^+_s,X^0_s,X^\pm_p\}$ and $Y_U$ the corresponding $\Uq$-module
(i.e. $R^+_s,U^+_s,$ $U^+_s,U^\pm_p$ respectively). Here $Y''=X^+_1$ unless $(Y,Y')=(X^0_s,P^+_\ell)$ or $(X^+_s,X^+_\ell)$, in which case $Y''=X^0_1$. Moreover,  $\gamma=\gamma(Y,Y')$ is given by the table
$$\vbox{\tabskip=0pt\offinterlineskip
  \def\tablerule{\noalign{\hrule}}
  \halign to 2in{
    \strut#&\vrule#\tabskip=0em plus1em &    
    \hfil#&\vrule#&\hfil#&\vrule#&\hfil#&\vrule#&     
\hfil#&\vrule#&\hfil#&
\vrule#\tabskip=0pt\cr
&&$Y\backslash Y'$&&\hfil$P^+_\ell$\hfil&&\hfil$X_\ell^+$ \hfil&&\hfil$X^0_\ell$\hfil&&\hfil$X^\pm_p$\hfil&\cr
\tablerule\tablerule&&$P^+_s\ $\hfil&&\hfil$\cP$\hfil&&\hfil$\cP$\hfil&&\hfil$\emptyset$\hfil&&\hfil$\cP$\hfil&\cr
\tablerule&&\hfil$X^+_s\ $\hfil&&\hfil$\cP$\hfil&&\hfil$\cP$\hfil&&\hfil$\emptyset$\hfil&&\hfil$\cP$\hfil&\cr
\tablerule&&\hfil$X^0_s\ $\hfil&&\hfil$x$\hfil&&\hfil$\emptyset$\hfil&&\hfil$\cP$\hfil&&\hfil$\emptyset$\hfil&\cr
\tablerule&&\hfil$X^\pm_p\ $\hfil&&\hfil$\cP$\hfil&&\hfil$\cP$\hfil&&\hfil$\emptyset$\hfil&&\hfil$\cP$
\hfil&\cr
\tablerule\noalign{\smallskip}}}$$
where as usual $1\le s,\ell<p$. (The choice of $\gamma$ is not always unique.)
For some examples, see \eqref{ex:S}; also, for the choice $(Y,Y')=(X^0_s,X^-_p)$ we obtain 
\[
\frac{S^\chi_{X^0_s,X^-_p}}{S^\chi_{X^+_1,X^-_p}}=\frac{S^{\hopflink}_{U^+_s,U^-_p}}{S^{\hopflink}_{U^+_1,U^-_p}}=0.
\]
We thus see that the logarithmic Hopf link invariants correspond to the $S$-matrix of an $\text{SL}(2, \mathbb Z)$-action, analogous to the Hopf link invariants in  modular tensor categories. 

In fact we can say more. Rescaling $x$ appropriately, it can be shown that there exist nonzero scalars $c^\chi(Y)$ and $c^\hopflink(\gamma)$
such that
\begin{equation}\label{cS=cS}c^\chi(Y') S^\chi_{Y,Y'}=c^\hopflink(\gamma)S^{\hopflink\,;\gamma}_{{Y_U},\bar{Y'_U}}\end{equation}
for all $Y,Y'$. This carries a little more information than \eqref{S/S=S/S}. At this point it isn't
clear to us whether the stronger and simpler \eqref{cS=cS} will continue to hold in other examples than $\cW_p$, or whether \eqref{S/S=S/S} is the  fundamental relation.

In Section 3.1, we gave different ways to obtain Verlinde-like formulas for 
strongly-finite VOAs. Let us work out explicitly the final possibility mentioned there, for $\cW_p$.
It gives a Verlinde-type formula for some of the tensor ring (as opposed to Grothendieck ring) structure constants. Namely, the $S$-matrix restricted to the $R_i$ is the $S$-matrix of theta functions as we will see below. It is in particular invertible and hence
\begin{equation}\label{eq:vertrip}
{\cN_{U^+_i, R_j}}^{R_k} = \sum_{\ell=0}^p \frac{S^{\hopflink\,;x}_{U^+_i, R_\ell}S^{\hopflink\,;\cP}_{R_j, R_\ell}  \left({S^{\hopflink\,;\cP}}^{-1}\right)_{R_\ell, R_k} }{S^{\hopflink\,;x}_{U^+_1, R_\ell}}.
\end{equation}

We conclude with the following observation from \cite{CG}:
\begin{theorem}\label{determined}
The tensor subring of $\Uq$, spanned by the simple modules $U^\pm_s$ and their projective covers
$R^\pm_s$, is completely determined from the following three datum: the socle series of projective modules;  the fact that the ambidextrous element $U^+_p$ squares to the tensor unit $U^+_1$; and the complete  list of logarithmic Hopf link invariants.
\end{theorem}

\subsection{Towards more examples}\label{sec:ex}

There is obviously a need for new examples, and  currently much effort is put into developing techniques to construct them \cite{CKL, CKLR, CKM}. Given a VOA one can construct a new one as a VOA extension, as a kernel of a screening charge, as an orbifold or a coset. We will state three constructions of these types that are expected to yield new examples of $C_2$-cofinite but non-rational VOAs.

\subsubsection{Higher rank analogues of $\cW_p$}

The triplet algebra $\cW_p$ is constructed as the kernel of a screening charge (the zero-mode of a specific intertwining operator) acting on the lattice VOA of the lattice $\sqrt{2p}\mathbb Z=\sqrt{p}A_1$. This has a fairly obvious generalization to root lattices of simply laced Lie algebras \cite{AM2, FT}. We follow \cite{CM2}. Let $Q$ be a root lattice of type $A, D$ or $E$ and let $L=\sqrt{p}Q$ for positive integer $p$ at least equal to two. Further denote by $\alpha_1, \dots, \alpha_n$ a set of simple roots of the Lie algebra $\mathfrak g$ corresponding to $Q$. Then one defines $\cW(p)_Q$ as the intersection of intertwiners associated to the simple roots:
\[
\cW(p)_Q:= \bigcap_{i=1}^n \text{Ker}\left( e_0^{-\frac{\alpha_i}{{\sqrt{p}}}} : V_L \rightarrow V_{L-\frac{\alpha_i}{\sqrt{p}}}\right).
\]
Characters of certain simple modules are known and they form a vector-valued modular form where the modular weights of homogeneous components ranges between zero and $|\Delta_+|$, the number of positive roots \cite{BM}. $C_2$-cofiniteness is not proven but believed to be true and the characters are of course an indication for this. 

Similar to the rank one case a close relation of the representation category of these VOAs and those of restricted quantum groups of type $\mathfrak g$ at $2p$-th root of unity is expected. As in the rank one case this won't be a braided equivalence. Since it seems much easier to study the quantum groups it is an instructive task to understand for example what happens in the case of $\mathfrak g=\mathfrak{sl}_3$ and compare it to $\cW(p)_{A_2}$.

\subsubsection{Logarithmic parafermion or $Z$-algebras}

Let $\cC_k$ be the  coset VOA $\cC_k=\text{Com}\left(\cH, L_k(\mathfrak{g})\right)$, where $L_k(\mathfrak{g})$ is the simple affine VOA of the simple Lie algebra $\mathfrak{g}$ at level $k$ and $\cH$ is the Heisenberg sub algebra of rank the rank of $\mathfrak g$. By a parafermion VOA we mean any VOA extension of $\cC_k$. 

Parafermion algebras \cite{FZ} have first been introduced by Lepowsky and Wilson under the name $Z$-algebras \cite{LW1}. They have been studied in the rational case, that is as Heisenberg cosets of rational affine VOAs. Under the assumption that the representation category of the parafermion VOA of non-integer admissible level is a vertex tensor category in the sense of \cite{HLZ8} a great deal of the representation theory of the parafermion algebra follows from the one of its affine parent theory \cite{CKLR}.

Let $\mathfrak g$ be a simple Lie algebra, then the WZW theory of $\mathfrak g$ at admissible rational negative level is surely not $C_2$-cofinite. However the space of highest weight vectors for the Heisenberg sub VOA corresponding to the Cartan subalgebra will be an abelian intertwining algebra whose fusion corresponds to addition in the root lattice of $\mathfrak g$. 
Restricting to a sublattice with the property that all involved Heisenberg highest-weight vectors have integral conformal dimension for the Virasoro field of the Heisenberg VOA gives then a VOA (under the vertex tensor category assumption). The outcome of \cite{CKLR} is that simple modules of this VOA are rare and we strongly believe that these VOAs furnish new examples of the finite logarithmic type. In the case of $\mathfrak g = \mathfrak{sl}_2$ this indeed works on the level of characters \cite{ACR}. There, we have no access to Hopf links so that its use as a further testing ground for our ideas is limited and will nonetheless be pursued.

\subsubsection{Non-trivial extensions of triplet algebras}

Another strategy to obtain new examples is as extensions of tensor products of triplet algebras, 
\[
\cW:=\cW_{p_1}\otimes \dots \otimes \cW_{p_n}.
\]
The precise question is: let $Y_1, \dots Y_s$ be a collection of simple (or even more general indecomposable) $\cW$-modules, then can the module
\[
Y_1 \oplus \dots\oplus Y_s
\]
be given the structure of a VOA? 

 The best accessible type are simple current extensions which indeed furnish new examples \cite{CKL}. Going beyond simple current extensions or abelian intertwining algebras is usually very difficult, but it can be done for some families \cite{CKM}. We also hope to test our ideas in that framework further. At the moment we only announce that for extensions by abelian intertwining algebras we have equally perfect agreement between Hopf links of the corresponding algebras obtained in the quantum group side and modular $S$-matrices on the VOA side.   
Finally, one can also look at extensions that are super VOAs, then modified super traces of open Hopf links will nicely agree with the modular properties of super characters and super pseudo-characters.

\subsection{Finite tensor categories and the Leavitt algebra}\label{sec:leavitt}

Recall that a modular tensor category is among other things a fusion category.
Inspired by the theory of subfactors, there is an interpretation for fusion categories using
algebra endomorphisms. You can think of it as a faithful representation of the category. 
The idea is that objects $M$ in the category (e.g. for our categories, the
$\cV$-modules) are realised by algebra homomorphisms $\phi_M:\cA\rightarrow\cA$
on some fixed algebra $\cA$. The point
is that tensor products $M\otimes N$ (the most difficult ingredient in the category) are now easy:
$\phi_{M\otimes N}=\phi_M\circ\phi_N$. This realisation by endomorphisms has evolved
into a powerful method for constructing and classifying fusion categories. 

These categories do not come naturally with a braiding. However, the centre- or double-construction
is a standard way to obtain a modular tensor category from a fusion category. For a category 
realised by endomorphisms, this construction reduces to linear algebra and is called the
tube algebra construction.

The endomorphism method has been used to construct several fusion categories, and 
through them several modular tensor categories, that seem much more difficult
to construct in other ways. It is known that any unitary fusion category can be realised by endomorphisms  on a $C^*$-algebra. Most modular tensor categories are not doubles
of fusion categories, though. The strongly-finite VOAs whose category of modules is a double of a fusion category,
are the strongly-rational VOAs which are subVOAs of holomorphic VOAs with the same central charge
(holomorphic VOAs are strongly-rational with a single irreducible module).

This endomorphism method came from subfactors, where it only worked for unitary fusion
categories. However, unitarity is a serious restriction, because semi-simplicity is then automatic. 
But \cite{EG} extended the method to nonunitary fusion categories, replacing $C^*$-algebras
with (uncompleted) Leavitt algebras. 
In \cite{CG} we explain how to extend \cite{EG} to realise finite tensor categories with
endomorphisms. Our hope is that this would be an effective way to generate new finite tensor
categories, which we then double to get new log-modular tensor categories. Just as one may 
expect that every modular tensor category is the category of modules of some strongly-rational
VOA, one may imagine that every log-modular tensor category  is the category of modules of some strongly-finite VOA. So we envisage this as a method for fishing for relatively simple strongly-finite
VOAs.

See \cite{CG} for a detailed nonsemi-simple example.

\subsection{Representation-type}\label{sec:reptype}

Let $\cV$ be strongly-finite and write $P=P_0\oplus\cdots\oplus P_n$
for the direct sum of projective covers of simple $\cV$-modules. Then $\cA_\cV:=\mathrm{End}_\cV(P)$ is a finite-dimensional
associative algebra, and $\mathrm{Mod}^{g.r}(\cV)$ is equivalent as an abelian
category to $\mathrm{Mod}^{fin}(\cA_\cV)$, where $\cV$-module $M$ corresponds to the right $\cA_\cV$-module $\mathrm{Hom}_\cV(P,M)$. Equivalence as abelian categories  means we ignore tensor products and duals, but the equivalence preserves simples,
projectives, indecomposables, composition series, dimensions of Hom-spaces, etc.

If we have in addition rigidity, then Mod$^{g.r}(\cV)$ is in fact a finite tensor category, and we can say more.
In particular, Proposition 2.7 of \cite{EO} says any finite tensor category can be regarded as the tensor category of modules of a
weak quasi-Hopf algebra. These are not as nice though as quasi-Hopf algebras: `weak' means
the counit is not a homomorphism of algebras and comultiplication is not unit-preserving.

\begin{proposition} {Suppose $\cV$ is strongly-finite and its category is rigid. Then ${\cV}$ has finitely many indecomposable
modules,  iff it is strongly-rational. If $\cV$ is
not strongly-rational, then $\cV$ has uncountably many indecomposable modules of arbitrarily high Jordan--H\"older
length.}\end{proposition}

The proof \cite{CG}  uses the Second Brauer--Thrall conjecture, proved by Nazarova--Roiter \cite{NR, ASS}.

We say something has \textit{tame} representation type if it  has finitely many 1-parameter families of indecomposable modules of each dimension.
 A great example of a tame-type algebra is the algebra $A=\bbC[x]$ of polynomials.
 Any $n$-dimensional module of $A$ is completely determined by how $x$ acts, which we can think of as an $n\times n$
 matrix $X$. Up to equivalence we can replace $X$ with a sum of Jordan blocks, as usual. An indecomposable module will then correspond to $X$ being a single Jordan block $X=B_{\lambda,n}$.
 Therefore for each dimension $n$, $X$ and hence the module is uniquely determined by the value $\lambda\in\bbC$,
 which can be arbitrary. This is the 1-parameter family of $n$-dimensional modules.
 Tame representation type isn't as easy as the finite-type representation theory of
 strongly-rational VOAs, but it is under control and all indecomposables could be classifiable.

 The remaining type is {\it wild-type}. This means that as the dimension $n$ grows, so does the number of parameters
 needed to parametrise the dimension $n$ indecoposables. The simplest and prototypical example is the group algebra $\bbC[F_2]$ of the free group $F_2=\langle
 a,b\rangle$. A representation corresponds to a choice of 2 invertible matrices $A,B$ of the same size $n\times n$:    
$a\mapsto A$ and $b\mapsto B$. Call this representation $\rho_{A,B}$.
Two representations $\rho_{A,B},\rho_{A',B'}$ are equivalent if there is an invertible matrix $P$ such that both
$A'=PAP^{-1}$ and $B'=PBP^{-1}$. It is easy to verify that there are $n^2+1$-dimensional families of pairwise inequivalent
$n$-dimensional indecomposable  $F_2$-representations. So it is indeed  much wilder than
the $\bbC[x]$ representation theory considered above. To parametrise all indecomposable representations of a wild-type algebra such as $\bbC[F_2]$ would be very complicated and probably useless.

The triplet algebras $\cW_p$ are tame. The symplectic fermions $SF^+_d$ (for $d>1$) are wild.
Almost all finite-dimensional associative algebras have wild representation type.
For these reasons, we'd expect almost all strongly-finite VOAs are wild type. 
Among other things, this means one should not try to classify their modules.

\section{The tensor category as a tool}\label{sec:tools}

A valid question is surely what the seemingly abstract nonsense about modular tensor categories is good for. As it furnishes a proof of the Verlinde formula in the rational case and the beginnings of a natural picture for it in the logarithmic case its use is already demonstrated. We now would like illustrate that the category also provides useful tools to answer questions about the structure of a VOA. 

\subsection{The VOA extension problem}

As discussed in the previous section, VOA extensions are natural tools to construct potential new $C_2$-cofinite but non-rational VOAs. 
The VOA extension problem has a concrete categorical counterpart, that of an algebra with unit. This is a recent development due to \cite{HKL}, see also \cite{KO}.

For this let $\mathcal C$ be a braided tensor category. An algebra $A$ in $\cC$ is 
an object $A$ together with multiplication $\mu$, that is a morphism
\[
\mu : \ A \otimes A \rightarrow A 
\]
and an embedding of tensor unit $\one$ of $\cC$, that is the VOA $\cV$ itself in the VOA setting:
\[
\iota: \ \one \hookrightarrow A.
\]
These have to be associative
\[
\mu \circ\left(\mu\otimes Id_A \right) \circ \cA = \mu \circ \left(Id_A \otimes \mu\right)\ : \ A \otimes \left(A\otimes A \right) \rightarrow A
\]
and commutative
\[
\mu \circ c_{A, A} =\mu
\]
and the unit condition 
\[
\mu\circ \left(\iota \otimes Id_A \right) \circ\ell_A^{-1}  = Id_A
\]
has to hold, with left-multiplication
\[
\ell_A : \ \one\times A \rightarrow A.
\]
Further the unit needs to be unique in the sense that Hom$_\cC\left(\one, A\right)$ is one-dimensional. 
The main theorem of \cite{HKL} is that VOA extensions of sufficiently nice VOAs are in one-to-one correspondence to algebras in the representation category of the VOA. 
There is one subtlety with the uniqueness of unit which is in practice fairly easy to check. Furthermore, the representation category of the extended VOA is isomorphic to the category of local algebra modules. As a consequence one can deduce a fair amount of general statements on the structure of VOA extensions as well as orbifolds and cosets. A few examples of interesting results that are proven employing VOA extensions are:
\begin{ex}
 Extensions by invertible objects (simple currents and more general abelian intertwining algebras) are still most accessible. Combining orbifold and coset they seem to give rise to new logarithmic $C_2$-cofinite VOAs. The picture \cite{CKL, CKLR} is as follows
\[
\cV \xrightarrow{\quad \text{Heisenberg}-\text{coset}\quad } \cC \xrightarrow{\quad \text{simple current extension}\quad } \cD. 
\]
Here $\cV$ is a fairly nice VOA containing a possibly higher rank Heisenberg VOA. Then under the assumption that $\cC$-modules are objects of a braided tensor category one can deduce that $\cD$ in the first place is a VOA and moreover its simple objects are classified by those of $\cV$. This is part of the outcome of \cite{CKLR} building on \cite{CKL}. The triplet algebras can all be realized in this picture starting with admissible level $W^{(2)}_{p-1}$-algebras for $\cV$ and the intermediate cosets $\cC$ are then the singlet algebras \cite{CRW}.

 If one takes for $\cV$ the VOA $L_k(\mathfrak{sl}_2)$ for negative admissible level $k$, that is $-2 < k < 0$ and $k$ rational and $k+2\neq 1/n$ for integer $n$, then $\cD$ has online finitely many simple objects and its character vector sublemented by pseudotrace functions is modular \cite{ACR}. But this strategy goes beyond cosets of affine VOAs and seems also to work for minimal W-algebras with one-dimensional associated variety \cite{ACKL}.
\end{ex}
\begin{ex}
Recently, Dong and Ren have proven rationality of parafermion VOAs of rational affine VOAs \cite{DoRe}, and this result has built on quite an amount of results \cite{ALY, DLWY}. 
The category picture provides a short and more general proof, that is given that $\cV$ is simple, rational and $C_2$-cofinite containing a lattice VOA of a positive definite lattice, then the coset is as well simple, rational and $C_2$-cofinite \cite{CKLR}. 
\end{ex}
\begin{ex}
A nice example that illustrates the power of the interplay of modularity and the tensor category is the recent classification of holomorphic VOAs of central charge $24$ \cite{vEMS}. There an orbifold construction of the remaining cases of Schelleken's famous list \cite{Sch} is given. The construction used that the orbifold VOA is rational \cite{CaMi} and hence its representation category is modular. Trace functions  of twisted modules are obtainable and their modular properties gave the Grothendieck ring of the orbifold VOA. The new holomorphic VOAs were then obtained as extensions of the orbifold VOAs by certain abelian intertwining algebras. 
\end{ex}
\begin{ex}
Let $\cV$ be a VOA, $\cW$ a subVOA and $\cC$ its commutant or coset VOA, such that $\cW $ and $\cC$ form a mutually commuting pair. 
Let 
\[
\cV = \bigoplus_{i\in I} \cW_i \otimes \cC_i
\]
the decomposition of $\cV$ as a $\cW\otimes \cC$-module. 
The mirror extension conjecture (respectively theorem in the regular case) states that if there is a subset $J\subset I$ and multiplicities $m_j$ such that
\[
\bigoplus_{j\in J} m_j \cW_j
\]
is a VOA, then so is
\[
\bigoplus_{j\in J} m_j \cC_j^\vee.
\]
The reason for this statement to be true in the regular setting, i.e. all three $\cW, \cC$ and $\cV$ are regular is, that the $\cW_i$ (respectively $\cC_i$) for $i$ in $I$ form a closed subcategory of the module category of $\cW$ (respectively $\cC$) and there is a braid-reversing equivalence between these two subcategories, relating $\cW_j$ to $\cC_j^\vee$. This is proven in \cite{Lin} and has been conjectured in \cite{DJX}.
In the non-rational setting a similar statement can be deduced if the fusion algebras of these two subcategories are both abelian intertwining algebras \cite{CKLR}.
\end{ex}

\subsection{Applications of Hopf links}

Obviously the main application of Hopf links is that they are the representation matrices of representations of the fusion ring. This fact has interesting consequences in studying VOAs used in \cite{CKL, CKLR}. Further the ideas below are important in understanding simple current extensions of tensor products of triplet algebras (and more generally abelian intertwining algebras). 

\subsubsection{Hopf links as parity indicators}

Let $\cV$ be a VOA and $J$ be a simple current of order two. Then a natural question is if 
\[
\cV\oplus J
\]
is a VOA or a super VOA or none of the two. Provided that the conformal dimension of $J$ is half-integer or integer one always obtains either a VOA or a super VOA \cite{CKL}. The parity question is decided by the quantum dimension of $J$,
\[
\text{qdim}(J)= \frac{S^\hopflink_{J, \cV}}{S^\hopflink_{\cV, \cV}} \ \in \ \{ \pm 1\}
\]
and one has the following cases
\begin{itemize}
\item $\cV\oplus J$ is a VOA if $\text{qdim}(J)=1$ and $J$ has integral conformal dimension;
\item $\cV\oplus J$ is a super VOA if $\text{qdim}(J)=1$ and $J$ has half-integral conformal dimension;
\item $\cV\oplus J$ is a super VOA if $\text{qdim}(J)=-1$ and $J$ has integral conformal dimension;
\item $\cV\oplus J$ is a VOA if $\text{qdim}(J)=-1$ and $J$ has half-integral conformal dimension;
\end{itemize}
While the first two cases are the standard VOA and super VOAs we would like to note that also the latter two cases are fairly rich. For example WZW theories of Lie superalgebras (affine super VOAs) are integer graded. A prominent example of the last type are the Bershadsky-Polyakov algebras \cite{Ber, Pol} and many other quantum Hamiltonian reductions. 

Parity of a simple current relates to their OPE (see Proposition 4.2 of \cite{CKL}) and this has very useful consequences. For example it turns out that this information toegther with the Jacobi identity uniquely specifies the OPE structure of a simple current extension of certain rational type $A$ W-algebras times lattice VOAs. As a consequence these extensions are rational Bershadsky-Polyakov algebras \cite{ACL}.

\subsubsection{Hopf links and locality} 

Let $\cV$ be a VOA and $J$ a simple current such that
\[
\cW= \bigoplus_{n} J^{\otimes n}
\]
is a VOA again. One then would like to understand which $\cV$-modules lift to modules of the extended VOA $\cW$. According to \cite{HKL} VOA extensions are in one-to-one correspondence with haploid algebras in the category and modules of the extension correspond to local modules of the algebra. 
Locality means that multiplication of the algebra on a module commutes with braiding. Again via the balancing axiom this can easily be related to properties of conformal weights: A module $X$ lifts to a local module 
\[
\mathcal F(X) = \bigoplus_n J^{\otimes n} \otimes X
\]
if and only if $h_{J\otimes X} -h_J- h_X$ is an integer. Here, also half-integral extensions are allowed and super VOAs as well. If $X$ is projective then this statement applies provided $X$ is a subquotient of the category generated by simple objects under tensor product.

\subsubsection{Hopf links and fixed-points of simple currents}

Let $J$ be a simple current of finite or infinite order. A module $X$ is called a fixed-point if $J^{\otimes s}\otimes X\cong X$ for some positive integer $s$ smaller than the order of $J$. Fixed-points of simple currents make the study of simple current extensions and orbifolds by lattice VOAs more complicated and one needs to be aware of then. They can be detected as follows:
The open Hopf link satisfies
\[
\Phi_{X, P} = \Phi_{J^{\otimes s}\otimes X, P} = \Phi_{J^{\otimes s}, P} \circ \Phi_{X, P}
\]
and thus either $\Phi_{X, P}=0$ or $\Phi_{J^{\otimes s}, P}=\left(\Phi_{J, P}\right)^s =\Id_P$. 
But using the balancing axiom it follows that
\[
t_P\left(\Phi_{J^{\otimes s}, P}\right) = \left(\text{qdim}(J)\right)^s d(P)\left(\theta_{J^{\otimes s} \otimes P} \circ \left(\theta^{-1}_{J^{\otimes s}} \circ \theta^{-1}_P\right)\right).
\]
Here $d(P)$ is the modified dimension of $P$, that is the modified trace of the identity on $P$. 
If on the other hand $\Phi_{X, P}=0$ then $0=t_P\left(\Phi_{X, P}\right) =S^{\hopflink\,;\cP}_{X, P}$. 
We thus get the following criterion: For any module $P$, $X$ being a fixed-point requires that at least one of the two holds:
\begin{itemize}
\item $S^{\hopflink\,;\cP}_{X, P}=0$;
\item $\left(\text{qdim}(J)\right)^s \left(\theta_{J^{\otimes s} \otimes P} \circ \left(\theta^{-1}_{J^{\otimes s}} \circ \theta^{-1}_P\right)\right)=1$,
\end{itemize}
Depending on $P$ this might sometimes be the first condition, sometimes the second and possibly even both.
An example of this is given in \cite{CKLR}.
The second condition is mainly a conformal dimension condition as the twist $\theta_X$ acts on $X$ as $e^{\pi i L_0}$.

\hspace*{1cm}

\noindent Department of Mathematical and Statistical Sciences, University of Alberta,
Edmonton, Alberta  T6G 2G1, Canada. 
\emph{email: creutzig@ualberta.ca and tjgannon@ualberta.ca}


\begin{thebibliography}{FGST2}








\bibitem[Ab]{Abe} T. Abe, \textit{A $\bbZ_2$-orbifold of the symplectic fermionic vertex operator superalgebra}, 
Math. Z. 255 (2007), 755--792.

\bibitem[AA]{AA} T. Abe and Y. Arike, \textit{Intertwining operators and fusion rules for vertex operator algebras arising from symplectic fermions,} J. Alg. 373 (2013), 39--64.

\bibitem[AdM1]{AM} D. Adamovic and A. Milas, \textit{On the triplet vertex algebra W(p)}, Adv. Math. 217 (2008) 2664--2699.

\bibitem[AdM2]{AM2} D. Adamovic and A. Milas, \textit{C2-cofinite vertex algebras and their logarithmic modules}, Proceedings of the conference Conformal Field Theories and Tensor Categories, Beijing Mathematical Lectures from Beijing University, Vol. 2, ( 2014), 18 pp.

\bibitem[AdM3]{AM3} D. Adamovic and A. Milas, \textit{An analogue of modular BPZ equation in logarithmic (super)conformal field theory}, Contemporary Mathematics (2009) 497, 1-17.

\bibitem[AdM4]{AM4} D. Adamovic and A. Milas, \textit{The structure of Zhu's algebras for certain W-algebras}, Advances in Mathematics , 227, (2011) 2425-2456.

\bibitem[AC]{AC} C. Alfes and T. Creutzig, \textit{The mock modular data of a family of superalgebras}, Proc. Amer. Math. Soc. 142 (2014), 2265-2280.

\bibitem[ACKL]{ACKL} T. Arakawa, T. Creutzig, K. Kawasetsu and A. Linshaw, \textit{Orbifolds and Cosets of minimal W-algebras}, in preparation.

\bibitem[ACL]{ACL} T. Arakawa, T. Creutzig, and A. Linshaw, \textit{Cosets of Bershadsky-Polyakov algebras and rational $\cW$-algebras of type $A$}, arXiv:1511.09143.

\bibitem[ACR]{ACR} J. Auger, T. Creutzig and D, Ridout, \textit{Modularity of logarithmic parafermion vertex algebras}, in preparation.

\bibitem[AnM]{AnMo} G. Anderson and G. Moore, \textit{Rationality in conformal field theory.} Commun. Math. Phys. 117 (1988), 441--450.

\bibitem[AN]{AN} Y. Arike and  K. Nagatomo, \textit{Some remarks on pseudo-trace functions for orbifold models associated with symplectic fermions.} Internat. J. Math. 24 (2013), no. 2, 1350008, 29 pp.

\bibitem[ALY]{ALY} T. Arakawa, C. H. Lam and H. Yamada, \textit{Zhu's algebra, C2-algebra and C2-cofiniteness of parafermion vertex operator algebras}, Adv. Math. 264 (2014), 261--295.

\bibitem[ASS]{ASS} I. Assem, D. Simson, and A. Skowronski, \textit{Elements of the representation theory of associative algebras, Vol. 1 and 2.} (Cambridge University Press, 2006).

\bibitem[Ber]{Ber} M. Bershadsky, \textit{Conformal field theories via Hamiltonian reduction},  139 (1991), no. 1, 71-82.

\bibitem[Bo]{Bo} R. E. Borcherds, \textit{Monstrous moonshine and monstrous Lie superalgebras}. Invent. Math. 109 (1992), no. 2, 405–-444.

\bibitem[BrV]{BrV} A. Brugui\`eres and A. Virelizer, \textit{Quantum double of Hopf monads and categorical centers}, Trans. Amer. Math. Soc. 365 (2012) 1225--1270.

\bibitem[BK]{BK} B. Bakalov and A. Kirillov, \textit{Lectures on tensor categories and modular functors}, University Lecture Series, 21. American Mathematical Society, (2001).

\bibitem[BM]{BM} K. Bringmann and A.Milas, \textit{W-Algebras, False Theta Functions and Quantum Modular Forms, I}, IMRN, 2015.

\bibitem[BW]{BW} J. Barrett and B. Westbury, \textit{Spherical categories}, Adv. Math. 143 (1999), 357--375.


\bibitem[CGP]{CGP} F. Costantino, N. Geer, and B. Patureau-Mirand, \textit{Some remarks on the unrolled quantum group of sl(2)}. J. Pure Appl. Algebra 219 (2015), no. 8, 3238--3262.

\bibitem[C]{C} T. Creutzig, \textit{Quantum dimensions in logarithmic CFT}, Oberwolfach report No. 16 (2015) 52--54.

\bibitem[CG]{CG} T. Creutzig and T. Gannon, \textit{The Theory of $C_2$-cofinite VOAs}, draft available at http://www.ualberta.ca/$\sim$creutzig. 

\bibitem[CKL]{CKL} T. Creutzig, S. Kanade and A. R. Linshaw, \emph{Simple current extensions beyond semi-simplicity}, \href{}{arXiv:1511.08754}, submitted.

\bibitem[CKM]{CKM} T. Creutzig, S. Kanade and R. McRae, \emph{Tensor categories for vertex operator algebra extensions}, in preparation.

\bibitem[CKLR]{CKLR} T. Creutzig, S. Kanade, A. R. Linshaw and D. Ridout, \emph{Schur-Weyl Duality for Heisenberg Cosets}, in preparation.

\bibitem[CL1]{CL} T. Creutzig and A. Linshaw, \textit{Cosets of affine vertex algebras inside larger structures},  arXiv:1407.8512. 

\bibitem[CL2]{CL2} T. Creutzig and A. Linshaw, \textit{The super $W_{1+\infty}$ algebra with integral central charge}, Trans. Amer. Math. Soc. 367 (2015), no. 8, 5521–-5551.

\bibitem[CL3]{CL3} T. Creutzig and A. Linshaw, \textit{Orbifolds of symplectic fermion algebras}, Trans. Amer. Math. Soc.,  arXiv:1404.2686.

\bibitem[CaMiy]{CaMi} S. Carnahan, M. Miyamoto, \textit{Regularity of fixed-point vertex operator subalgebras}, arXiv:1603.05645.

\bibitem[CM1]{CM1} T. Creutzig and A. Milas, \textit{False theta functions and the Verlinde formula}, Adv. Math. {\bf 262} (2014) 520--545.

\bibitem[CM2]{CM2} T. Creutzig and A. Milas,  \textit{Higher rank partial and false theta functions and representation theory}, in preparation.

\bibitem[CMR]{CMR} T. Creutzig, A. Milas and M. Rupert, \textit{Logarithmic Link Invariants of $\overline{U}^H_q(\mathfrak{sl}_2)$ and Asymptotic Dimensions of Singlet Vertex Algebras}. in preparation.

\bibitem[CMW]{CMW} T. Creutzig, A. Milas and S. Wood, \textit{On regularised quantum dimensions of the singlet vertex operator algebra and false theta functions}, arXiv:1411.3282.

\bibitem[CQS1]{CQS1} T. Creutzig, T. Quella and V. Schomerus, \textit{T. Creutzig, T. Quella and V. Schomerus, Branes in the GL(1$|$1) WZNW-Model}, Nucl. Phys. B792 (2008) 257-283.

\bibitem[CQS2]{CQS} T. Creutzig, T. Quella and V. Schomerus, \textit{New boundary conditions for the c=-2 ghost system}, Phys. Rev. D77 (2008) 026003.

\bibitem[CRo]{CRo} T. Creutzig and P. B. Ronne, \textit{The GL(1$|$1)-symplectic fermion correspondence}, Nucl. Phys. B815 (2009) 95-124.

\bibitem[CR1]{CR1} T. Creutzig and D. Ridout, \textit{Modular data and Verlinde formulae for fractional level WZW models I}, Nucl. Phys. B {\bf 865} (2012) 83--114.

\bibitem[CR2]{CR2} T. Creutzig and D. Ridout, \textit{Modular data and Verlinde formulae for fractional level WZW models II}, Nucl. Phys. B {\bf 875} (2013) 423--458.

\bibitem[CR3]{CR3} T. Creutzig and D. Ridout, \textit{Relating the archetypes of logarithmic conformal field theory}, Nucl. Phys. B {\bf 872} (2013) 348-391.

\bibitem[CR4]{CR4} T. Creutzig and D. Ridout, \textit{Logarithmic conformal field theory: beyond an introduction}, J.  Phys. A {\bf 46} (2013) 4006.

\bibitem[CRW]{CRW} T. Creutzig, D. Ridout and S. Wood, \textit{Coset Constructions of Logarithmic (1,p)-Models}, Lett. Math. Phys. 104, 5 (2014) 553-583.

\bibitem[DJX]{DJX} C. Dong, X. Jiao and F. Xu, \textit{Mirror extensions of vertex operator algebras}, Comm. Math. Phys. 329 (2014), 263-294.


\bibitem[DS]{DS} B. J. Day and R. Street, \textit{Centres of monoidal categories of functors}, Contemp. Math. 431 (2007) 187--202.

\bibitem[DoRe]{DoRe} C. Dong and L. Ren, \textit{Representations of the parafermion vertex operator algebras}, arxiv:1411.6085.

\bibitem[DR]{DR} A.~Davydov and I.~Runkel, \textit{$\bbZ/2\bbZ$-extensions of Hopf algebra module categories by their base categories}, Adv. Math. 247 (2013) 192--265. 




\bibitem[DLWY]{DLWY} C. Dong, C. H. Lam. Q. Wang and H. Yamada, \textit{The structure of parafermion vertex operator algebras}, J. Algebra 323 (2010), no. 2, 371–-381. 

\bibitem[EGNO]{EGNO} P. Etingof, S. Gelaki, D. Nikshych, and V. Ostrik, \textit{Tensor categories} (American Math. Soc.,
Providence 2015).

\bibitem[vEMS]{vEMS} J. van Ekeren, S. M\"oller, N. R. Scheithauer, \textit{Construction and Classification of Holomorphic Vertex Operator Algebras}, arXiv:1507.08142.

\bibitem[EG1]{EG0} D.E. Evans and T. Gannon, \textit{The exoticness and realisability
of twisted Haagerup-Izumi modular data}, \textit{Commun. Math. Phys.}  307 (2011), 463--512. 

\bibitem[EG2]{EG} D.E. Evans and T. Gannon, \textit{Non-unitary fusion categories and their doubles via endomorphisms}, arXiv:1506.03546.

\bibitem[ENO]{ENO}  P. Etingof, D. Nikshych and V. Ostrik, \textit{An analogue of Radford's S4 formula for finite tensor categories}. Int. Math. Res. Not. 2004, no. 54, 2915–-2933.

\bibitem[EO]{EO} P. Etingof and V. Ostrik, \textit{Finite tensor categories},  Mosc. Math. J. {\bf 4} (2004) 627--654, 782--783. 

\bibitem[F]{F} E. Frenkel, \textit{Langlands correspondence for loop groups}. Cambridge Studies in Advanced Mathematics, 103. Cambridge University Press, Cambridge, 2007.

\bibitem[FGST1]{FGST1} B.~L.~Feigin, A.~M.~Gainutdinov, A.~M.~Semikhatov and I.~Y.~Tipunin, \textit{Modular group representations and fusion in logarithmic conformal field theories and in the quantum group center}, Commun. Math. Phys.  {\bf 265} (2006) 47--93.




\bibitem[FGST2]{FGST2} B.~L.~Feigin, A.~M.~Gainutdinov, A.~M.~Semikhatov and I.~Y.~Tipunin, \textit{The Kazhdan-Lusztig correspondence for the representation category of the triplet W-algebra in logarithmic CFT}, 
  Theor. Math. Phys. {\bf 148} (2006) 1210--1235.



\bibitem[Fi]{Fi} F. Fiordalisi, \textit{Logarithmic intertwining operators and genus-one correlation functions}, arXiv:/1602.03250.

\bibitem[FHL]{FHL} I. B. Frenkel, Y.-Z.Huang and J. Lepowsky, \textit{On axiomatic approaches to vertex operator algebras and modules}, Memoirs Amer. Math. Soc. 104, 1993.

\bibitem[FHST]{FHST}  J. Fuchs, S. Hwang, A. M. Semikhatov, and I. Y. Tipunin,  \textit{Nonsemisimple 
fusion algebras and the Verlinde formula}, Commun. Math. Phys. 247 (2004) 713--742.

\bibitem[FK]{FK} M. Flohr, H. Knuth, \textit{On Verlinde-Like Formulas in c(p,1) Logarithmic Conformal Field Theories}, arXiv:0705.0545.

\bibitem[FRS]{FRS}  J. Fuchs, I. Runkel and C. Schweigert, \textit{TFT construction of RCFT correlators. I -- IV}, Part I Nuclear Phys. B 646 (2002), no. 3, 353--497; Part II, Nuclear Phys. B 678 (2004), no. 3, 511-- 637; Part III, Nuclear Phys. B 694 (2004), no. 3, 277--353; Part IV, Nuclear Phys. B 715 (2005), no. 3, 539--638.

\bibitem[FFRS]{FFRS}  J. Fuchs, J. Fjelstad, I. Runkel and C. Schweigert, \textit{TFT construction of RCFT correlators. V. Proof of modular invariance and factorization}. Theory Appl. Categ. 16 (2006), No. 16, 342-–433.

\bibitem[FLM]{FLM}  I. Frenkel, J. Lepowsky and A. Meurman, \textit{Vertex operator algebras and the Monster}. Pure and Applied Mathematics, 134. Academic Press, Inc., Boston, MA, 1988.  

\bibitem[FS1]{FS} J. Fuchs and C. Schweigert, \textit{Hopf algebras and finite tensor categories in conformal field
theory.}  {Rev. Un. Mat. Argentina} 51 (2010), 43--90.

\bibitem[FS2]{FS2} J. Fuchs and C. Schweigert, \textit{Consistent systems of correlators in non-semisimple conformal field theory}, arXiv1604.01143. 

\bibitem[FS3]{FS3} J. Fuchs and C. Schweigert, \textit{Coends in conformal field theory}, arXiv:1604.01670.

\bibitem[FSS]{FSS} J. Fuchs, C. Schweigert and C. Stigner \textit{From non-semisimple Hopf algebras to correlation functions for logarithmic CFT}, J. Phys. A: Math. Theor. 46 (2013) 494008.

\bibitem[FT]{FT} B. Feigin and I. Tipunin, \textit{Logarithmic CFTs connected with simple Lie algebras}, arXiv:1002.5047.

\bibitem[FZ]{FZ} V. A. Fateev and A. B. Zamolodchikov, \textit{Nonlocal (parafermion) currents in two-dimensional conformal quantum
field theory and self-dual critical points in ZN -symmetric statistical systems}, Sov. Phys. JETP 62, 215--225 (1985).

\bibitem[GabG]{GG} M. R. Gaberdiel and	P. Goddard, \textit{Axiomatic conformal field theory}, Commun.Math.Phys. 209 (2000) 549--594.

\bibitem[GabK]{GK} M. R. Gaberdiel and H. G. Kausch, \textit{A rational logarithmic conformal field theory}, Phys.Lett. B386 (1996) 131--137.

\bibitem[GabR]{GR} M. R. Gaberdiel and I. Runkel, \textit{From boundary to bulk in logarithmic CFT}, J. Phys. A 41 (2008) 075402.


\bibitem[GaiR1]{GaR} A.~M.~Gainutdinov and I.~Runkel, \textit{Symplectic fermions and a quasi-Hopf algebra structure on $\bar{U}_i sl(2)$}, arXiv:1503.07695.

\bibitem[GaiR2]{GaR2} A.~M.~Gainutdinov and I.~Runkel, \textit{The non-semisimple Verlinde formula and pseudo-trace functions} arXiv:1605.????

\bibitem[GaiT]{GaiT} A. M. Gainutdinov and I.Y. Tipunin, \textit{Radford, Drinfeld, and Cardy boundary states in (1,p) logarithmic conformal field models}, J. Phys. A 42 (2009) 315207. [0711.3430[hep-th]].

\bibitem[GKP1]{GKP1} N. Geer, J. Kujawa, and B. Patureau-Mirand, \textit{Ambidextrous objects and trace functions for nonsemisimple categories}. Proc. Amer. Math. Soc. 141 (2013), no. 9, 2963--2978.

\bibitem[GKP2]{GKP2} N. Geer, J. Kujawa, and B. Patureau-Mirand, \textit{Generalized trace and modified dimension functions on ribbon categories}. Selecta Math. (N.S.) 17 (2011), no. 2, 453--504.

\bibitem[GPT]{GPT} N. Geer, B. Patureau-Mirand, and V. Turaev, \textit{Modified quantum dimensions and re-normalized link invariants}. Compos. Math. 145 (2009), no. 1, 196--212.

\bibitem[GPV]{GPV} N. Geer, B. Patureau-Mirand, and A. Virelizier, \textit{Traces on ideals in pivotal categories}. Quantum Topol. 4 (2013), no. 1, 91--124.


\bibitem[H1]{H1} Y.-Z. Huang, \textit{Vertex operator algebras and the Verlinde conjecture}, Commun. Contemp. Math. 10 (2008) 103--154.

\bibitem[H2]{H2} Y.-Z. Huang, \textit{Rigidity and modularity of vertex tensor categories}, Commun. Contemp. Math. 10 (2008) 871--911. 


\bibitem[HKL]{HKL} Y.-Z. Huang, A. Kirillov Jr. and J. Lepowsky, \textit{Braided Tensor Categories and Extensions of Vertex Operator Algebras}, Comm. Math. Phys. 337 (2015), no. 3, 1143--1159.



\bibitem[HL]{HL3} Y.-Z. Huang and J. Lepowsky, \textit{A theory of tensor products for module categories for a vertex operator algebra. IV}. J. Pure Appl. Algebra 100 (1995) 1-3, 173--216.
			








\bibitem[HLZ]{HLZ8} Y.-Z. Huang, J. Lepowsky, and L. Zhang, \textit{Logarithmic tensor category theory, I -- VIII}, arXiv:1012.4193v7, arXiv:1012.4196v2, arXiv:1012.4197v2, arXiv:1012.4198v2,
arXiv:1012.4199v3, arXiv:1012.4202v3, arXiv:1110.1929v3, arXiv:1110.1931v2.

\bibitem[JS]{JS} A. Joyal, R. Street, \textit{Braided Tensor Categories}, Advances in Mathematics, 102, 1, (1993), 20--78.

\bibitem[Kas]{K} C. Kassel, \textit{Quantum groups}. Graduate Texts in Mathematics, 155. Springer-Verlag, New York, 1995.

\bibitem[Kau]{Kau} H. G. Kausch, \emph{Symplectic fermions.} {Nucl.\ Phys.} {B583} (2000), 513--541.

\bibitem[KL]{KL} T. Kerler and V.V. Lyubashenko, \textit{Non-semisimple topological quantum field theories for 3-manifolds with corners.} Springer Lecture Notes in Mathematics 1765 (Springer Verlag, New York 2001)

\bibitem[KO]{KO} A.\ Kirillov Jr., V. Ostrik, \emph{On a $q$-analogue of
  the McKay correspondence and the ADE classification of
  $\mathfrak{sl}_2$ conformal field theories}. {Adv.\ Math.} {
    171} (2002), no.\ 2, 183--227.

\bibitem[KS]{KS} H. Kondo and Y. Saito, \textit{Indecomposable decomposition of tensor products of modules over the restricted quantum universal enveloping algebra associated to sl2}. J. Algebra 330 (2011), 103--129.

\bibitem[Lin]{Lin}X. Lin, \textit{Mirror extensions of rational vertex operator algebras}, arXiv:1411.6742.

\bibitem[LM]{LM} R. Laber and G. Mason, \textit{$\bbC$-graded vertex algebras and conformal flow}, arXiv:1308.0557.

\bibitem[LW]{LW1} J. Lepowsky and R. L. Wilson, \textit{A new family of algebras underlying the Rogers-Ramanujan identities and generalizations}, Proc. Natl. Acad. Sci. USA 78, 7254--7258 (1981).




\bibitem[Ly1]{L} V. Lyubashenko, \textit{Modular transformations for tensor categories}. J. Pure Appl. Algebra 98 (1995), no. 3, 279--327. 

\bibitem[Ly2]{L2} V. Lyubashenko, \textit{Modular properties of ribbon abelian categories}, arXiv:hep-th/9405168.

\bibitem[Ly3]{L3} V. Lyubashenko, \textit{Invariants of 3-manifolds and projective representations of mapping class
groups via quantum groups at roots of unity}, Commun.Math. Phys. 172 (1995) 467–-516.


\bibitem[Mil1]{Mil} A. Milas, \textit{Weak modules and logarithmic intertwining operators for vertex operator algebras},
In: Recent developments in infinite-dimensional Lie algebras and conformal field theory (Charlottesville, VA, 2000), Contemp. Math., 297, (Amer. Math. Soc., Providence, RI, 2002) pp. 201--225.

\bibitem[Mil2]{Mil2} A. Milas, \textit{Characters of modules of irrational vertex operator algebras }, Conformal Field Theory, Automorphic Forms and Related Topics CFT 2011, Heidelberg, Vol.8 (2014), 1-17.

\bibitem[Miy1]{M0} M. Miyamoto, \textit{Intertwining operators and modular invariance.} Surikaisekikenkyusho Kokyuroku No. 1218 (2001), 57–-68.


\bibitem[Miy2]{M1} M. Miyamoto, \textit{Modular invariance of vertex operator algebra satisfying $C_2$-cofiniteness}, Duke Math. J. 122 (2004) 51--91.

\bibitem[Miy3]{M2} M. Miyamoto, \textit{A theory of tensor products for vertex operator algebra satisfying $C_2$-cofiniteness}, arXiv:0309350v3

\bibitem[Miy4]{M3} M. Miyamoto, \textit{Flatness and semi-rigidity of vertex operator algebras}, arXiv:1104.4675.

\bibitem[Miy5]{Miy4} M. Miyamoto, \textit{C2 -Cofiniteness of cyclic-orbifold models}, Commun. Math. Phys. {\bf 335} (2015) 1279--1286.

\bibitem[MS]{MS} G. W. Moore and N. Seiberg, \textit{Classical and quantum conformal field theory}, Commun. Math. Phys. {\bf 123} (1989) 177--254.

\bibitem[MSV]{MSV} F. Malikov, V. Schechtman and A.Vaintrob, \textit{Chiral de Rham complex}. Comm. Math. Phys. 204 (1999), no. 2, 439-–473.

\bibitem[Mug]{Mue} M. M\"uger, \textit{Tensor categories: a selective tour}, Rev. Un. Mat. Argent. {\bf 51}
 (2010)  95--163.


\bibitem[Mur]{Mu} J. Murakami, \textit{From colored Jones invariants to logarithmic invariants}. arXiv:1406.1287.

\bibitem[MN]{MN} J. Murakami and K. Nagatomo, \textit{Logarithmic knot invariants arising from restricted quantum groups}. Internat. J. Math. 19 (2008), no. 10, 1203--1213.

\bibitem[NR]{NR} L. A. Nazarova and A. V. Roiter, \textit{Kategorielle Matrizen-Probleme und die
Brauer-Thrall-Vermutung}. Mitt. Math. Sem. Giessen Heft 115 (1975).

\bibitem[NT]{NT} K. Nagatomo and A. Tsuchiya, \textit{The triplet vertex operator algebra W(p) and 
the restricted quantum group at root of unity}, Adv. Stud. in Pure Math., Amer. Math. Soc. 61 (2011) 1--49.

\bibitem[Pol]{Pol} A. Polyakov, \textit{Gauge transformations and diffeomorphisms}, Internat. J. Modern Phys. A 5 (1990), no. 5, 833-842.

\bibitem[PRR]{PRR} P. A. Pearce, J. Rasmussen, and P. Ruelle, \textit{Grothendieck ring and Verlinde-like formula for the $\cW$-extended logarithmic minimal model $\mathcal{WLM}(1,p)$}, J. Phys. A43 (2010) 045211, 13pp. 

\bibitem[R]{R} J. Rasmussen, \textit{Fusion matrices, generalized Verlinde formulas, and partition functions in $\mathcal{WLM}(1,p)$},  J. Phys. A 43 (2010), no. 10, 105201, 27 pp.

\bibitem[RW1]{RW1} D.~Ridout and S.~Wood, \textit{Modular transformations and Verlinde formulae for logarithmic $(p_+,p_-)$-models}, Nucl. Phys. B {\bf 880} (2014) 175--202.

\bibitem[RW2]{RW2} D.~Ridout and S.~Wood, \textit{The Verlinde formula in logarithmic CFT}, J. Phys. Conf. Ser.  {\bf 597} (2015) 1,  012065.

\bibitem[Ru]{Ru} I. Runkel, \textit{A braided monoidal category for free super-bosons,} J. Math. Phys. 55 (2014), no. 4, 041702, 59 pp.

\bibitem[Sch]{Sch} A. N. Schellekens, \textit{Meromorphic c = 24 conformal field theories}, Comm. Math. Phys. 153 (1993), 159–-185.

\bibitem[Sh]{Sh} K. Shimizu, \textit{The monoidal center and the character algebra,} arXiv:1504.01178v2.

\bibitem[Su]{S} R. Suter, \textit{Modules over $U_q(sl_2)$}. Comm. Math. Phys. 163 (1994), no. 2, 359--393.

\bibitem[TW1]{TW} A. Tsuchiya and S. Wood, \textit{The tensor structure on the representation category of the $\cW_p$ triplet algebra},  J. Phys. A 46 (2013), no. 44, 445203, 40 pp.

\bibitem[TW2]{TW2} A. Tsuchiya and S. Wood, \textit{On the extended W-algebra of type $sl_2$ at positive rational level}, International Mathematics Research Notices, 2014


\bibitem[T]{T} V. G.Turaev, \textit{Quantum invariants of knots and 3-manifolds}. de Gruyter Studies in Mathematics, 18. Walter de Gruyter $\&$ Co., Berlin, 1994. x+588 pp.

\bibitem[V]{V}   E. P. Verlinde, \textit{Fusion rules and modular transformations in 2D conformal field theory}, Nucl. Phys. B {\bf 300} (1988) 360--376.

\bibitem[Z]{Z} Y. Zhu, \textit{Modular invariance of characters of vertex operator algebras}, J. Amer. Math. Soc. 9 (1996) 237--302.

\end{thebibliography}
\end{document}